\input ./style/arxiv-general.cfg
\documentclass[MSNbibl,number,citesort,dvips]{arxbj}
\makeatletter
   \@ifpackageloaded{graphicx}{}{\usepackage{graphicx}}
\makeatother
\usepackage{mathbh}

\aid{0}
\volume{22}
\issue{1}
\pubyear{2016}
\firstpage{396}
\lastpage{420}
\doi{10.3150/14-BEJ663} 
\docsubty{FLA}

\makeatletter
\newcommand{\iint}{\int\!\!\!\int}
\newcommand{\rrvert}{\vert}
\newcommand{\rrVert}{\Vert}
\newcommand{\llvert}{\vert}
\newcommand{\llVert}{\Vert}
\newcommand{\eqref}[1]{(\ref{#1})}
\newtheorem{theorem}{Theorem}
\newtheorem{lemma}{Lemma}
\newremark{remark}{Remark}
\newcommand{\nn}{\nonumber}
\newcommand{\Ind}{\mathbh{1}}
\renewcommand{\P}{\mathrm{P}}
\makeatother

\begin{document}
\begin{frontmatter}

\title{Adaptive Bayesian density regression for high-dimensional data}
\runtitle{Bayesian density regression in high dimensional}

\begin{aug}
\author[1]{\inits{W.}\fnms{Weining}~\snm{Shen}\corref{}\thanksref{1}\ead[label=e1]{wshen@mdanderson.org}}
\and
\author[2]{\inits{S.}\fnms{Subhashis}~\snm{Ghosal}\thanksref{2}\ead[label=e2]{sghosal@ncsu.edu}}
\address[1]{Department of Biostatistics, The University of Texas MD
Anderson Cancer Center, 1400 Pressler Street, Houston, TX
77030, USA. \printead{e1}}
\address[2]{Department of Statistics, North Carolina State University,
4276 SAS Hall, 2311 Stinson Drive, Raleigh, NC 27695, USA. \printead{e2}}
\end{aug}

\received{\smonth{7} \syear{2013}}
\revised{\smonth{6} \syear{2014}}

\begin{abstract}
Density regression provides a flexible strategy for modeling the
distribution of a response variable $Y$ given predictors $\mathbf{X} =
(X_1,\ldots,X_p)$ by letting that the conditional density of $Y$ given
$\mathbf{X}$ as a completely unknown function and allowing its shape to
change with the value of $\mathbf{X}$. The number of predictors $p$ may be
very large, possibly much larger than the number of observations $n$,
but the conditional density is assumed to depend only on a much smaller
number of predictors, which are unknown. In addition to estimation, the
goal is also to select the important predictors which actually affect
the true conditional density. We consider a nonparametric Bayesian
approach to density regression by constructing a random series prior
based on tensor products of spline functions. The proposed prior also
incorporates the issue of variable selection. We show that the
posterior distribution of the conditional density contracts adaptively
at the truth nearly at the optimal oracle rate, determined by the
unknown sparsity and smoothness levels, even in the ultra
high-dimensional settings where $p$ increases exponentially with $n$.
The result is also extended to the anisotropic case where the degree of
smoothness can vary in different directions, and both random and
deterministic predictors are considered. We also propose a technique to
calculate posterior moments of the conditional density function without
requiring Markov chain Monte Carlo methods.
\end{abstract}

\end{frontmatter}

\section{Introduction}\label{sec1}
We consider Bayesian estimation of the conditional density of a
response $Y$ given a large number of predictors $\mathbf{X}=(X_1,\ldots
,X_p)$, where $p$ is possibly much larger than the sample size $n$.
This problem is sometimes referred to as density regression and has
received attention in many scientific application areas such as genome
association studies. Non-Bayesian approaches to density regression
usually focus on the kernel approach \cite{Fan1996,Hall1999}, which
requires estimating the bandwidth using cross-validation \cite{Fan2004},
bootstrap \cite{Hall1999} or other methods.

In the Bayesian literature, there are two common approaches to density
regression. One approach models the joint density and obtains the
conditional density as a by-product. The other approach directly models
the conditional density while leaving the marginal distribution of $\mathbf
{X}$ unspecified. In this paper, we focus on the latter approach in a
nonparametric manner. Many of the existing methods are based on
assigning priors on the space\vadjust{\goodbreak} of densities through countable mixtures
of location-scale densities \cite{Muller1996} or through generalizing
stick-breaking representations \cite{chungdunson09,dunsonpark08,Dunson2007,Griffin2006}. Priors obtained by
transforming a Gaussian process \cite{Jara2011,Tokdar2010} and a
multivariate generalization of a beta process \cite{Trippa2011} have
also been considered. Ma \cite{Ma2012} proposed a generalized P\'{o}lya
tree, which possesses nice posterior conjugacy properties and hence
allows fast computation.

In modern data analysis, often the data may be high-dimensional.
Statistical analysis in such a setting is possible only under some
sparsity assumption and only if a variable selection procedure is
implemented. Many variable selection techniques have been introduced in
the frequentist literature, such as discrete subset selection and
penalization methods. Popular methods include the least absolute
shrinkage and selection operator (lasso) introduced in \cite{Tib1996}
and the sure independence screening (SIS) proposed in \cite{Fan2008}.
Under the $p \gg n$ setting, oracle properties of lasso-type estimators
have been established for parametric models including linear regression
in \cite{Greenshtein2004}, generalized linear model in \cite{vandeGeer2008} and for nonparametric additive models in \cite{Huang2010}. For nonparametric (conditional) density estimation
problems, however, similar results are only obtained under a fixed $p$
setting in \cite{Hall2004,Lepski2012}.



Bayesian variable selection methods have also gained popularity.
For example, stochastic search variable selection (SSVS) adopts an
efficient sampling-based method to avoid comparing all possible
sub-models \cite{Carlin1995,George1993,George1997}.
Bayesian model averaging methods incorporate model uncertainty into
estimation and predictions \cite{Brown1998,Brown2002}.
Bayesian variable selection is commonly accomplished by assigning a
Bernoulli distribution prior on each covariate \cite{chungdunson09},
whereas an efficient SSVS algorithm is implemented to search the model
space and to combine the posterior estimation results from different
models. Tokdar \textit{et al.} \cite{Tokdar2010} extended variable selection to
dimension reduction by allowing the true sets of covariates determined
by a sub-linear space of $\mathbf{X}$ through a projection operator. While
these proposed methods show promising numerical results, rates of
contraction are largely unknown. Moreover, modern applications often
require that we allow the dimension $p$ of the predictor to be also
large, possibly much larger than $n$. So far such results are largely
missing from the Bayesian nonparametric literature.

In the linear regression problem, recovery of the regression
coefficients requires nontrivial assumptions on the structure of
covariates, for example, the restricted isometry property or the
compatibility condition to make the underlying problem well posed; see
\cite{Buhlmann2011}, although the corresponding problem of estimating
the regression function does not require such conditions, see, for
example, \cite{Dalalyan2012} for a discussion under a Bayesian
framework. In the density regression context, the recovery of the
conditional density is analogous to that of the regression function in
the linear regression context and hence does not require such conditions.

In the recent years, the literature on Bayesian asymptotics has
flourished with many fundamental breakthroughs. General results for
posterior contraction rates were established in \cite{Ghosal2000,Ghosal20071,Ghosal2007,Shen2001,Vander2008,Rivoirard2012}.
For density regression models, a consistency result was obtained by
Tokdar \textit{et al.} \cite{Tokdar2010} for a logistic Gaussian process prior,
by Norets and Pelenis \cite{Norets2012} for a kernel stick-breaking
process prior and by Pati \textit{et al.} \cite{Pati2011} for a probit
stick-breaking process prior. Tokdar \cite{Tokdar2011} also obtained
the posterior convergence rate for the logistic Gaussian process prior
given a fixed $p$. For high-dimensional Bayesian models, there are very
few contraction rates results available. Parametric models have been
studied by Jiang \cite{Jiang2007} for generalized linear model and by
Castillo and van der Vaart \cite{Castillo20121} for Gaussian white
noise model. A classification model with categorical predictors was
considered by Yang and Dunson \cite{Yang2012}, who constructed priors
using tensor factorizations and obtained a posterior contraction rate
allowing $p$ to grow exponentially with $n$.

In this paper, we consider the Bayesian density regression problem
using a finite linear combination of tensor products of B-splines to
construct a prior distribution. We obtain the posterior contraction
rate under the $p \gg n$ setting and show that the rate is adaptive for
both dimension and smoothness in the sense that it agrees with the
optimal rate of convergence (up to a logarithmic factor) of the oracle
procedure, that uses the knowledge of true predictors and the
underlying smoothness of the true conditional density, simultaneously
for all smoothness levels and dimension of true predictors. We further
extend the result to the anisotropic situation where smoothness can
vary in different coordinates of the true conditional density function,
and allow both random and deterministic predictors.
We also devise an effective computing strategy. Typically, a reversible
jump Markov chain Monte Carlo (RJMCMC) introduced by \cite{Green1995}
is used for Bayesian computation for models with varying dimension. For
high dimensional data, RJMCMC-based methods may be computationally
challenging and may give unreliable results due to limited exploration
of the model space. We propose a group-wise Dirichlet distribution
prior on the coefficients of B-spline functions that leads to a
conjugacy-like structure which can be utilized to develop a computing
algorithm based on direct sampling without resorting to MCMC
techniques. As in the univariate density estimation example in \cite{Shen2012}, the proposed computing method presents closed form
expressions for posterior moments including the mean and the variance.

The paper is organized as follows. In the next section, we describe the
model and the prior and discuss some preliminaries on tensor product of
B-splines. Posterior contraction rates for both the isotropic and
anisotropic cases and for both random and deterministic predictors are
obtained in Section~\ref{secrate}. Computational strategies are described and
simulation results are presented in Section~\ref{secnum}. Proofs are presented in
Section~\ref{sec5}.

\section{Bayesian density regression}\label{sec2}

\subsection{Notation}\label{sec21}

Let $\mathbb{N} = \{1,2,\ldots\}$, $\mathbb{N}_0 = \{0,1,2,\ldots\}$
and $\Delta_J$ be the unit $J$-dimensional simplex. For any real number
$x$, define $\lfloor x \rfloor$ to be the largest integer less than or
equal to $x$. For a multi-index $\mathbf{l}=(l_1,\ldots,l_d) \in\mathbb
{N}_0^d$, $d \in\mathbb{N}$, we define the sum $\mathbf{l}.=l_1+\cdots
+l_d$ and the mixed partial derivative operator $D^{\mathbf{l}}=\partial
^{\mathbf{l}.}/\partial x_1^{l_1}\cdots\,\partial x_d^{l_d}$. For a bounded
open connected set $\Omega \subset\mathbb{R}^d$ (e.g., $\Omega=
(0,1)^d$), define the $\beta$-H\"{o}lder class $\mathcal{C}^{\beta
}(\Omega)$ as the collection of functions $f$ on $\Omega$ that has
bounded mixed partial derivatives $D^{\mathbf{l}} f$ of all orders up to
$\mathbf{l}. \leq\beta_0$ and that for every $\mathbf{l} \in\mathbb{N}_0^d$
satisfying $\mathbf{l}. = \beta_0$,
%
\begin{equation}
\bigl|D^{\mathbf{l}}f(\mathbf{x}) - D^{\mathbf{l}}f(\mathbf{y}) \bigr| \leq C \|
\mathbf{x}-\mathbf{y}\| _2^{\beta- \beta_0}
\end{equation}
for some constant $C>0$, any $\mathbf{x}, \mathbf{y} \in\Omega$ and $\beta_0$
as the largest integer strictly smaller than $\beta$. Any such function
uniquely extends to a continuous function on the closure of $\Omega$.





Let the indicator function be denoted by $\Ind$. We use ``$\lesssim$''
to denote an inequality up to a constant multiple. We write $f \asymp
g$ if $f \lesssim g \lesssim f$.
Let $D(\varepsilon,T,\rho)$ denote the packing number, which is defined as
the maximum cardinality of an $\varepsilon$-dispersed subset of $T$ with
respect to distance $\rho$.
The symbol $\P$ will stand for a generic probability measure.


\subsection{B-spline and its tensor-products}\label{sec22}

B-spline functions and their tensor-products have been widely used to
approximate functions in both mathematics and statistics literature.
Here we provide a brief overview of their definitions and approximation
properties; see more descriptions in \cite{Deboor2001}. For natural
numbers $K$ and $q \in\mathbb{N}$, let the unit interval $(0,1)$ be
divided into $K$ equally spaced subintervals. A spline of order~$q$
with knots at the end points of these intervals is a function $f$ such
that the restriction of~$f$ in each subinterval is a polynomial of
degree less than $q$ and if $q\ge2$, $f$ is $(q-2)$-times continuously
differentiable (interpreted as only continuous if $q=2$). Splines of
order $q$ form a linear space of dimension $J=K+q-1$, a convenient
basis of which is given by the set of B-splines $B_1,\ldots,B_J$. In
particular, if $q=1$, then corresponding B-splines form the Haar basis
$\{\Ind_{(j-1/J,j/J]}\dvt 1\le j\le J\}$. The B-splines are nonnegative and
add up to one at any $x$. Each $B_j$ is positive only on an interval of
length $q/K$ and at most $q$ many B-splines are nonzero at any given~$x$. Most importantly, splines of order $q$ with knots at $\{
0,1/K,\ldots,(K-1)/K,1\}$ approximate any function in $C^\alpha(0,1)$
at the rate $K^{-\alpha}$, or equivalently, any $f\in C^\alpha(0,1)$
can be approximated by a linear combination of $B_1,\ldots,B_J$ up to
an error of the order $J^{-\alpha}$.

This idea works in multidimensional case as well. For $(0,1)^d$ and $d
\in\mathbb{N}$, we split $(0,1)$ into $K_i$ equal intervals and
consider corresponding spline functions $B_1,\ldots,B_{{J_i}}$ in the
$i$th direction, where $J_i = q + K_i -1$. Hence, there are $\prod_{i=1}^d K_i$ equal cubes in total. Define tensor-product B-spline
basis functions as the product of univariate basis functions of each direction:
%
\begin{equation}
B_{\mathbf{j}} (x_1,\ldots,x_d) = \prod
_{k=1}^d B_{j_k} (x_k),
\qquad \mathbf{j}=(j_1,\ldots,j_d), j_k=1,
\ldots,J_k, k=1,\ldots,d.
\end{equation}
For simplicity, we use $\mathbf{B}$ to denote a column vector of all basis
functions and define the total number of basis functions by $J =\prod_{k=1}^d J_k$.

Tensor-products of B-splines maintain a lot of nice properties that the
univariate B-splines enjoy. In the following, we list a few of them
that will be used in our modeling:
\begin{itemize}[(iii)]
\item[(i)] $0 \leq B_{\mathbf{j}} \leq1$, for every $\mathbf{j} = (j_1,\ldots
,j_d) \in\{1,\ldots,J_1\} \times\cdots\times\{1,\ldots,J_d\}$.
\item[(ii)] $\sum_{j_1=1}^{J_1} \cdots\sum_{j_d=1}^{J_d} B_{\mathbf{j}}(\mathbf
{x}) =1$, for\vspace*{1.5pt} every $\mathbf{x} \in(0,1)^d$.
\item[(iii)] For every $\mathbf{x}$, $B_{\mathbf{j}}(\mathbf{x}) > 0$ only if
$\lfloor x_i K \rfloor\leq j_i \leq\lfloor x_i K\rfloor+ q-1$ for
every $i=1,\ldots,d$.
\end{itemize}

We also define the normalized version of a univariate B-spline $B$ by $
\bar{B} = B/\int_0^1 B(x)\,\mathrm{d}x$. Like univariate B-splines, the
approximation ability of tensor-product B-splines is determined by the
smoothness level $\alpha$ of the function to be approximated and $J$
provided that $q$ is chosen to be larger than $\alpha$. In the
following lemma, we state their approximation results. In particular,
the result in part (c) suggests that the approximation power remains
the same when the coefficients satisfy certain restrictions (positive,
adds up to one), which later can help us assign prior distributions.

\begin{lemma}\label{lemmal00}
\textup{(a)} For any function $f \in\mathcal{C}^{\beta} ((0,1)^{d} )$,
$0 < \beta\leq q$, there exists $\bolds{\theta}\in\mathbb{R}^{J}$ and a
constant $C_1 > 0$ such that
\[
\Biggl\|f - \sum_{j_1=1}^{J_0} \cdots\sum
_{j_{d}=1}^{J_0} \theta_{\mathbf{j}} \mathbf{B_j}(\mathbf{x}) \Biggr\|_{\infty} \leq C_1
J_0^{-\beta} \bigl\|f^{(\beta)}\bigr\|_{\infty},
\]
where $\mathbf{j}=(j_1,\ldots,j_{d})$.

\textup{(b)} Further, if $f > 0$, then for sufficiently large $J_0$, we can
choose every element of $\bolds{\theta}$ to be positive.

\textup{(c)} Assume that $f(y|x_1,\ldots,x_d)$ is a positive density function\vspace*{1pt} in
$y$ for every $(x_1,\ldots,x_d)$ and as a function of $(y,x_1,\ldots
,x_d)$ belongs to $\mathcal{C}^{\beta} ((0,1)^{d+1} )$, where
$0 < \beta\leq q$. Then for sufficiently large $J = J_0^{d+1}$, there
exists $\bolds{\eta} \in(0,1)^{J}$ satisfying $\sum_{j_0=1}^{J_0} \eta
_{j_0,j_1,\ldots,j_d} =1$ for every fixed $(j_1,\ldots,j_d) \in\{
1,\ldots,J_0\}^d$, and a constant $C_2 > 0$ such that
\[
\Biggl\|f(y|x_1,\ldots,x_d) - \sum
_{j_0=1}^{J_0}\cdots\sum_{j_{d}=1}^{J_0}
\eta_{j_0,\ldots,j_d} \bar{B}_{j_0}(y) \prod
_{k=1}^{d} B_{j_k}(x_k)\Biggr\|
_{\infty} \leq C_2 J_0^{-\beta}=C_2
J^{-\beta/(d+1)}.
\]
\end{lemma}

\subsection{Data generating process and the prior}\label{sec23}

We consider the data generated from $n$ independent and identically
distributed pairs of observations $(Y_1,\mathbf{X}_1),\ldots,(Y_n,\mathbf
{X}_n)$, where $Y_i \in(0,1)$ and $\mathbf{X}_i \in(0,1)^p$ for every
$i=1,\ldots,n$, and $p \in\mathbb{N}$. It may be noted that the unit
intervals appearing in the ranges of these random variables are not
special as we can apply an appropriate affine transform on the data
otherwise. We assume that $Y$ is related only to $d$ covariates, say
$X_{m_1^0},\ldots,X_{m_d^0}$, that is the conditional density of~$Y$
given $\mathbf{X}$ is a function of these coordinates only. This is an
important sparsity assumption that will allow us to make valid
inference about the conditional density even when $p$ is very large,
provided that $d$ is small. However, neither $d$ nor these indexes are
known. The goal is to estimate the conditional density of $Y$ given $\mathbf
{X}$ with accuracy comparable with the oracle procedure which assumes
the knowledge of $d$ and $m_1^0,\ldots, m_d^0$ using a Bayesian procedure.

A prior on the conditional density given $p$ covariate values
$x_1,\ldots,x_{p}$ can be induced by
a finite series expansion in terms of tensor product of B-splines
%
\begin{equation}
\label{eqe100}
h(y,\mathbf{x}|J_0,\mathbf{J},\bolds{\eta}) = \sum
_{j_0=1}^{J_0}\cdots\sum
_{j_{p}=1}^{J_{p}} \eta_{j_0,\ldots,j_p}
\bar{B}_{j_0}(y) \prod_{k=1}^{p}
B_{j_k}(x_k),
\end{equation}
where $\mathbf{j}=(j_1,\ldots,j_{p})$ is a $p$-dimensional index,
$\mathbf{J}=(J_1,\ldots,J_{p})$ and $\bolds{\eta}=(\eta_{1,\mathbf{j}},\ldots,\eta
_{J_0,\mathbf{j}})^T$ lies in a $J_0$-dimensional simplex
for every $\mathbf{j}\le\mathbf{J}$.
Note that $\eta_{j_0,\ldots,j_p}$ does not change for $j_k=1,\ldots,J_k$,
if and only if the $k$th component does not affect the conditional
density. In order to incorporate this feature in the prior, we define
variable inclusion indicators
$\gamma_k = \Ind$ (the $k$th variable is in the model). Let $\bolds{\gamma}=(\gamma_1,\ldots,\gamma_p)$.
Thus, $\eta_{j_0,\ldots,j_p}$ depends only on $j_0$ and $j_k$ with $k\in
\operatorname{Supp}(\bolds{\gamma})
=\{k\dvt  \gamma_k=1\}=\{m_1,\ldots,m_r\}$ for $r=\sum_{k=1}^p \gamma_k$.
Thus, the common value of $\eta_{j_0,\ldots,j_p}$ can be denoted by
$\theta_{j_0,j_{m_1},\ldots,j_{m_r}}\!$.
Now the conditional density can be written as
%
\begin{equation}
\label{eqe1000}
h(y,\mathbf{x}|J_0,\mathbf{J},\bolds{\eta}) = \sum
_{r=0}^p \sum
_{m_1,\ldots,m_r} \sum_{j_0=1}^{J_0}
\sum_{j_{m_1}=1}^{J_{m_1}}\cdots\sum
_{j_{m_r}=1}^{J_{m_r}} \theta_{j_0,j_{m_1},\ldots,j_{m_r}} \bar
{B}_{j_0}(y) \prod_{k\dvt  \gamma_k=1}
B_{j_k}(x_k).
\end{equation}
By assigning prior distributions on indicator variables $\gamma_k$,
number of terms $J_k$ and the corresponding coefficients $\bolds{\theta}$,
we obtain an induced prior on $f$.
The prior on the model indicator $\bolds{\gamma}$ is constructed by first
putting a prior on the total model size $r$, and then selecting models
with size~$r$. More specifically, we construct the prior distribution
through the following scheme:

%
\begin{enumerate}[(A3)]
%
\item[(A1)] Prior on the model size $r$: Let $r= \sum_{k=1}^p \gamma_k$
be the number of variables included in the model and $\Pi_1$ be a
fixed, positive prior probability mass function of $r$. Assume that
there exists some constants $c_0, t_0 > 0$, such that for every $r \in
\mathbb{N}$,
%
\begin{equation}
\label{a1} \Pi_1(r) \leq\exp\bigl\{- \exp\bigl(c_0
r^{t_0}\bigr) \bigr\}.
\end{equation}
\item[(A2)] Prior on the inclusion variables $(\gamma_1,\ldots,\gamma
_p)$: Given a value of $r$, define the support of~$\bolds{\gamma}$ by $\{
m_1,\ldots,m_r\}$. We assume that the probability $\Pi_2(m_1,\ldots
,m_r|r)$ of each set of variables $1\le m_1<\cdots<m_r \le p$ of size
$r$ satisfies
\[
c_1' \frac{1}{{p \choose r}} \leq \Pi_2(m_1,
\ldots,m_r|r) \leq c_1''
\frac{1}{{p \choose r}}
\]
for some positive constant $c_1' \le c_1''$.
\item[(A3)] Prior on the number of terms in the basis expansion: Given
the model size $r$ and active predictor indices $1 \leq m_1 < \cdots<
m_r \leq p$, let the number of terms in the basis expansion in
$(r+1)$-fold tensor products of B-splines be denoted by
$(J_0,J_{m_1},\ldots,J_{m_r})$, and let $\Pi_3(\cdot|r;m_1,\ldots,m_r)$
stand for their joint prior distribution. We let $\Pi_3$ be induced by
independently distributed $J_0,J_{m_1},\ldots,J_{m_r}$
with identical distribution $\tilde{\Pi}_3$, that is,
\[
\Pi_3(J_0,J_{m_1},\ldots,J_{m_r} |
r;m_1,\ldots,m_r) = \tilde{\Pi }_3(J_0|r)
\prod_{k=1}^r \tilde{
\Pi}_3(J_{m_k}|r)
\]
and that for some fixed constants $c_2',c_2'' > 0, \kappa' \geq\kappa
'' \geq1$,
%
\begin{eqnarray}
\label{A3} \exp\bigl\{ -c_2' j^{r+1} (\log
j)^{\kappa'} \bigr\} \leq\tilde{\Pi}_3(j|r) \leq \exp\bigl\{
-c_2'' j^{r+1} (\log
j)^{\kappa''} \bigr\}.
\end{eqnarray}
%
%
\item[(A4)] Prior on the coefficients: Given the values of $r$,
$m_1,\ldots,m_r$ and $J_0,J_{m_1},\ldots,J_{m_r}$, recall that the
conditional density of $Y$ given $\mathbf{X}$ is written as
\begin{eqnarray*}
&& h(\mathbf{x},y |r;m_1,\ldots,m_r;J_0,J_{m_1},
\ldots,J_{m_r};\bolds {\theta})
\\
&& \quad = \sum_{j_0=1}^{J_0} \sum
_{j_{m_1}=1}^{J_{m_1}} \cdots \sum_{j_{m_r}=1}^{J_{m_r}}
\theta_{j_0,j_{m_1},\ldots,j_{m_r}} \bar {B}_{j_0}(y) B_{j_{m_1}}(x_{m_1})
\cdots B_{j_{m_r}}(x_{m_r}),
\end{eqnarray*}
where $(\theta_{j_0,j_{m_1},\ldots,j_{m_r}}\dvt 1 \leq j_0 \leq J_0) \in
\Delta_{J_0}$ for every $j_{m_1},\ldots,j_{m_r}$. We let every $(\theta
_{j_0,j_{m_1},\ldots,j_{m_r}}\dvt 1 \leq j_0 \leq J_0) \in\Delta_{J_0}$
be distributed independently with identical prior distribution $\tilde
{\Pi}_4(\cdot|J_0)$ and then denote the induced prior on the
coefficients by $\Pi_4(\cdot|r; m_1,\allowbreak\ldots,m_r; J_0,J_{m_1},\ldots
,J_{m_r})$. In particular, we choose $\tilde{\Pi}_4$ to be a Dirichlet
distribution $\operatorname{Dir}(a,\ldots,a)$, where $a$ is a fixed positive constant.
\end{enumerate}
Our prior (A1) on the model size includes the truncated binomial prior
used in \cite{Jiang2007} as a special case. Condition \eqref{a1}
implies that $\Pi_1(r > \bar{r}) \leq\exp\{-\exp(c_0' \bar{r}^{t_0})\}
$ for some constant $c_0' > 0$ and any $\bar{r} \in\mathbb{N}$. Since
$r$ should not be greater than $p$, which changes with $n$, we may also
let $\Pi_1$ depend on~$n$. In this case, we assume the decay rate to
hold with $c_0$ and $t_0$, which are both free from $n$, and for any
fixed $\bar{d}$, $\Pi_1(\bar{d})$ is bounded below, or more generally
$\Pi_1(\bar{d})$ satisfies $-\log\Pi_d(\bar{d}) = \mathrm{o}(n^{\delta})$ for
all $\delta> 0$.

In (A2), an easy choice is to let $c_1' = c_1''=1$, i.e., assign equal
probability for choosing $r$ indices from $\{1,\ldots,p\}$. We may also
allow $c_1',c_1''$ depend on $n$ as long as $\log(1/c_1') = \mathrm{o}(n^{\delta
})$ and $\log(c_1'') = \mathrm{o}(n^{\delta})$ for all $\delta> 0$. The
posterior contraction rate to be obtained in the next section will
remain the same.

For (A3), using the same prior $\tilde{\Pi}_3$ is not necessary, but it
is a convenient and appropriate default choice. The independence
between components is also not essential. In Section~\ref{anrate} when
the true density function is anisotropic, we shall have to use a
different and more complicated prior, which will be obviously
appropriate for this isotropic case as well. Relation \eqref{A3} is
satisfied with $\kappa'=\kappa''=1$ if a zero-truncated Poisson
distribution is assigned on $J^{r+1}$ in the sense that $N$ is a
zero-truncated Poisson and $J=\lfloor N^{1/(r+1)}\rfloor$.

In (A4), the same value of $a$ is not necessary, but we use it as a
default choice. In particular, $a=1$ leads to the uniform prior
distribution on the simplex. The same contraction rate will be obtained
as long as the parameters are kept in a fixed compact subset of
$(0,\infty)$. More generally, we may allow the lower bound approach
zero at most polynomially fast in $n^{-1}$, although the upper bound
needs to be fixed.


\section{Posterior contraction rates}\label{secrate}
\subsection{Isotropic case}\label{sec31}

In this section, we establish results on posterior contraction rates
for density regression. We allow the total number of covariates $p$
diverge with the increasing sample size $n$.
Let $\Pi$ be the prior as defined in (A1)--(A4) and denote the
posterior distribution based on $n$ pairs of observations $(Y_1,\mathbf
{X}_1), \ldots, (Y_n,\mathbf{X}_n)$ by $\Pi(\cdot| \mathbf{X}^n,\mathbf{Y}^n)$.
Let $\varepsilon_n\to0$ be a sequence of positive numbers. Consider a
suitable metric on the space of probability densities on $(0,1)$, such
as the Hellinger metric. Let $G$ stand for the common distribution of
$\mathbf{X}_1,\ldots,\mathbf{X}_n$, which need not be known.
We define the root average squared Hellinger distance on the space of
conditional densities by
%
\begin{equation}
\rho^2(f_1,f_2) = \iint \bigl
\{f_1^{1/2}(y|x_1,\ldots,x_p) -
f_2^{1/2}(y|x_1,\ldots,x_p) \bigr
\}^2 \,\mathrm{d}y G(\mathrm{d}x_1,\ldots,\mathrm{d}x_p),
\end{equation}
where $f_1$ and $f_2$ stand for generic conditional densities of $Y$ on
$(0,1)$ given $\mathbf{X}$ in $(0,1)^p$.
Let $f_0$ be a fixed conditional density function for $Y$ on $(0,1)$
given $\mathbf{X}$ in $(0,1)^p$, standing for the true conditional density.
We say that the posterior distribution of the density regression model
based on $\Pi$ contracts to $f_0$ at a rate $\varepsilon_n$ in the metric
$\rho$ if for any $M_n \rightarrow\infty$,
%
\begin{equation}
\label{defrate}
\lim_{n \to\infty} \Pi \bigl[\bigl\{f\dvt
\rho(f_0, f) > M_n \varepsilon_n\bigr\} | \mathbf
{X}^n, \mathbf{Y}^n \bigr] = 0\qquad \mbox{in probability}.
\end{equation}

We make the following assumptions.


\begin{enumerate}[(B3)]
\item[(B1)] The true density $f_0$ depends only on $d$ predictors
$X_{m_1^0},\ldots,X_{m_d^0}$, where $d$ is a fixed number. Further we
assume that as a function of $y$ and $x_{m_1^0},\ldots,x_{m_d^0}$, we
have $f_0 \in\mathcal{C}^{\beta} ((0,1)^{d+1} )$ for some $ 0
< \beta\leq q$.
\item[(B2)] The ambient dimension $p \leq\exp(C n^{\alpha})$ for $0 <
\alpha< 1$. When $\alpha= 0$, we interpret this condition as $p \leq
n^{K}$ for some constant $K > 0$.
\item[(B3)] The true conditional density $f_0$ is bounded below by a
positive constant $\underline{m}_0$.
\end{enumerate}

\begin{theorem}\label{thm}
Suppose that we have i.i.d. observations $\mathbf{X}_1,\ldots,\mathbf{X}_n$
from a possibly unknown probability distribution $G$ on $(0,1)^{p}$.
Assume that the true conditional density satisfies conditions
\textup{(B1)}--\textup{(B3)}. If the prior satisfies conditions \textup{(A1)}--\textup{(A4)}, then the
posterior distribution of $f$ contracts at $f_0$ at the rate
%
\begin{equation}
\label{isotropicrate}
\varepsilon_n = \max\bigl\{n^{-(1-\alpha)/2} (\log
n)^{1/(2t_0)}, n^{-\beta
/(2\beta+d+1)} (\log n)^{\kappa' \beta/(2\beta+ d+1)}\bigr\}
\end{equation}
with respect to $\rho$, where $t_0$ and $\kappa'$ are defined in \textup{(A1)}
and \textup{(A3)}.
\end{theorem}

Theorem~\ref{thm} establishes $\varepsilon_n$ as a bound on the posterior
contraction rate at $f_0$. It is known that the minimax rate associated
with estimating a $(d+1)$-dimensional density lying in a $\beta$-H\"
older class is $(n/\log n)^{-\beta/ (2\beta+ d +1)}$ (see \cite{Hasminskii1978}) with respect to the supremum norm. The minimax rate
of convergence for conditional densities with respect to the metric
$\rho$ is not known yet, but it is reasonable to expect that the rate
$n^{-\beta/ (2\beta+ d +1)}$ up to a logarithmic factor applies in
this situation as well, and can be taken as the oracle rate with which
the rates obtained in Theorem~\ref{thm} can be compared. Thus
if $p$ grows polynomially fast in $n$, then the rate we obtained
coincides with the oracle rate up to a logarithmic factor. If $p$ grows
exponentially fast, then it makes an impact on the rate. Note that we
obtain the optimal rate with the use of the same prior distribution for
all values of $\alpha$ and $\beta$. Hence our estimation and variable
selection procedure is rate-adaptive in the sense that the posterior
automatically adapts to the unknown number of covariates $d$ (i.e., the
oracle dimension) in the true model and the smoothness level $\beta$.
Our result also trivially contains the fixed dimensional situation
where no variable selection is involved. Note that the contraction at
the true density does not necessarily guarantee the convergence of the
selected set of predictors to the true set of predictors. The question
of recovering the true set of predictors remains open and is beyond the
scope of the present paper. However, as contraction rates are regulated
by the complexity of the underlying model determined by its dimension,
it may be anticipated that the posterior distribution assigns most of
its mass to low complexity models relative to the ambient dimension.

\begin{remark}
Theorem~\ref{thm} establishes contraction rates for the posterior
distribution of the entire conditional density function $f(y|\mathbf{x})$.
As a consequence, we can obtain the same posterior contraction rate for
other quantities of interest such as conditional quantile functions,
conditional moment functions and so on. This rate may not be optimal
for the estimation of these\vspace*{1pt} quantities because $y$\vspace*{1pt} has been integrated
out, that is, we conjecture the optimal rate is $n^{-\beta/(2\beta+
d)}$ instead of $n^{-\beta/(2\beta+ d +1)}$, up to logarithmic factors.
\end{remark}

\begin{remark}
After examining the proof, we find that condition \eqref{a1} in \textup{(A1)}
can be relaxed if $\alpha$ is small. For example, if $\alpha= 0$, then
we only need $\Pi_1(r) \leq\exp(-c_0 r^{t_0})$.
\end{remark}

\subsection{Anisotropic case}\label{anrate}
If predictors are qualitatively different, then it may be interesting
to consider the situation where $f_0$ has different smoothness levels
in different directions. In the following, we propose an alternative
anisotropic smoothness assumption replacing condition (B1).

For $\bolds{\beta} = (\beta_0,\ldots,\beta_d) \in\mathbb{N}^{d+1}$ and
$\beta_0,\ldots,\beta_d \leq q$, define a tensor Sobolev space $\mathcal
{S}^{\bolds{\beta}} ((0,1)^{d+1} )$ of functions $f$ of $(d+1)$
variables by
\[
\mathcal{S}^{\bolds{\beta}} \bigl((0,1)^{d+1} \bigr)= \bigl\{f\dvt
\bigl\|D^{\mathbf{l}} f \bigr\|_{\infty} < \infty, \mathbf{l}=(l_0,
\ldots,l_d), l_k \leq\beta_k, k=0,\ldots,d
\bigr\}
\]
with an associated norm $\|\cdot\|_{\mathcal{S}^{\bolds{\beta}}}$ defined as
\[
\|f\|_{\mathcal{S}^{\bolds{\beta}}} = \|f\|_{\infty} + \biggl\llVert \frac
{\mathrm{d}^{\beta_{0}} f}{\mathrm{d}y^{\beta_{0}}}
\biggr\rrVert _{\infty} + \sum_{k=1}^d
\biggl\llVert \frac{\mathrm{d}^{\beta_k} f}{\mathrm{d}x_k^{\beta_k}} \biggr\rrVert _{\infty}.
\]

As in Lemma~\ref{lemmal00}, we show that the tensor-product B-splines
still have nice approximation abilities within anisotropic function spaces.

\begin{lemma}\label{lemmal2}
\textup{(a)} For any function $f \in\mathcal{S}^{\bolds{\beta}}
((0,1)^{d+1} )$, where $0 < \beta_0,\ldots,\beta_d \leq q$, there
exists $\bolds{\theta}\in\mathbb{R}^{\prod_{k=0}^{d} J_k}$ and a constant
$C_1 > 0$ depending only on $q$, $d$ and $\beta_0,\ldots,\beta_d$ such that
\[
\Biggl\|f - \sum_{j_0=1}^{J_0} \cdots\sum
_{j_{d}=1}^{J_d} \theta_{\mathbf{j}} \mathbf
{B_j}(\mathbf{x}) \Biggr\|_{\infty} \leq C_1 \sum
_{k=0}^d J_k^{-\beta_k}
\biggl\llVert \frac{\mathrm{d}^{\beta_k} f}{\mathrm{d}x_k^{\beta_k}} \biggr\rrVert _{\infty},
\]
where $\mathbf{j}=(j_0,\ldots,j_{d})$.

\textup{(b)} Further, if $f > 0$, we can choose every element of $\bolds{\theta}$
to be positive.

\textup{(c)} Assume that $f(y|x_1,\ldots,x_d)$ is a positive probability density
function in $y$ for every $(x_1,\ldots,x_d)$ and as a function of
$(y,x_1,\ldots,x_d)$ belongs to $\mathcal{S}^{\bolds{\beta}}
((0,1)^{d+1} )$, where $\bolds{\beta}=(\beta_0,\ldots,\beta_d) \in
\mathbb{N}^{d+1}$ satisfying $0 < \beta_0,\ldots,\beta_{d} \leq q$.
Then there exists $\bolds{\eta} \in(0,1)^{\prod_{k=0}^d J_k}$ satisfying
$\sum_{j_0=1}^{J_0} \eta_{j_0,j_1,\ldots,j_d} =1$ for every fixed
$(j_1,\ldots,j_d) \in\{1,\ldots,J_1\} \times\cdots\times\{1,\ldots
,J_d\}$ and a~constant $C_2 > 0$ such that
\[
\Biggl\|f(y|x_1,\ldots,x_d) - \sum
_{j_0=1}^{J_0}\cdots\sum_{j_d=1}^{J_0}
\eta _{j_0,\ldots,j_d} \bar{B}_{j_0}(y) \prod
_{k=1}^{d} B_{j_k}(x_k) \Biggr\|
_{\infty} \leq C_2 \sum_{k=0}^d
J_k^{-\beta_k}.
\]
\end{lemma}

\begin{enumerate}[(B4)]
\item[(B4)] We assume that the true density $f_0$ is only related to
$d$ predictors with indices $1\le{m_1^0}<\cdots<{m_d^0}\le p$, where
$d$ is a fixed number, and as a function of $(y,x_{m_1^0},\ldots
,x_{m_d^0} )$ belongs to $\mathcal{S}^{\bolds{\beta}}
((0,1)^{d+1} )$.
\end{enumerate}
%
In order to obtain the adaptive convergence rate, we replace the
independent prior distribution on~$J$ in condition (A3) by the
following joint distribution condition.
\begin{enumerate}[(A3$'$)]
\item[(A3$'$)] Prior on the number of terms in basis expansion: Given
the model size $r$ and active predictor indices $1 \leq m_1 < \cdots<
m_r \leq p$, let the number of terms in the basis expansion of
$(r+1)$-fold tensor products of B-splines be denoted by
$(J_0,J_{m_1},\ldots,J_{m_r})$, and let $\Pi_3(\cdot|r;m_1,\ldots,m_r)$
stand for their joint prior distribution. We assume that for some fixed
constants $c_2',c_2'' > 0, \kappa' \geq\kappa'' \geq1$,
\begin{eqnarray*}
&& \exp \Biggl\{ -c_2' J_0 \prod
_{k=1}^r J_{m_k} \Biggl(\log
J_0 + \sum_{k=1}^r \log
J_{m_k} \Biggr)^{\kappa'} \Biggr\}
\\
&&\quad  \leq \Pi_3(J_0,J_{m_1},
\ldots,J_{m_r}|r;m_1,\ldots,m_r)
\\
&&\quad  \leq \exp \Biggl\{ -c_2''
J_0 \prod_{k=1}^r
J_{m_k} \Biggl(\log J_0 + \sum
_{k=1}^r \log J_{m_k}
\Biggr)^{\kappa''} \Biggr\}.
\end{eqnarray*}
\end{enumerate}
Then we obtain the posterior convergence rate for anisotropic functions.

\begin{theorem}\label{thm2}
Suppose that we have i.i.d. observations $\mathbf{X}_1,\ldots,\mathbf{X}_n$
from an unknown probability distribution $G$ on $(0,1)^{p}$. Assume
that the true conditional density satisfies conditions \textup{(B2)}--\textup{(B4)}. If
the prior satisfies conditions \textup{(A1)}, \textup{(A2)}, \textup{(A3$'$)} and \textup{(A4)}, then the
posterior distribution of $f$ contracts at $f_0$ at the rate
%
\begin{equation}
\label{anisotropicrate}
\varepsilon_n = \max \bigl\{n^{-(1-\alpha)/2} (\log
n)^{1/(2t_0)}, n^{-\beta^*/(2\beta^* + d + 1)} (\log n)^{\kappa' \beta^*/(2\beta^* +
d + 1)} \bigr\}
\end{equation}
with respect to $\rho$, where $\beta^* = (d+1) (\sum_{k=0}^d \beta
_k^{-1})^{-1}$ is the harmonic mean of $\{\beta_0,\ldots,\beta_d\}$,
$t_0$ and $\kappa'$ are defined in \textup{(A1)} and \textup{(A3$'$)}.
\end{theorem}

Clearly, the rate reduces to that of the isotropic case when $\beta
_0=\cdots=\beta_d$.
Thus the rate now can be viewed as the optimal rate (up to a
logarithmic factor) of estimating a $(d+1)$-variate function with
smoothness levels $\beta_0,\ldots,\beta_d$. Note that the rate is
determined by the harmonic mean of smoothness levels in different
coordinates, thus the worst smoothness has the most prominent effect on
the rate. However, the rate thus obtained is strictly better than that
obtained by a naive application of Theorem~\ref{thm} using the
worst smoothness condition in all co-ordinates. Thus additional
smoothness in other co-ordinates help improve the rate from the worst
case scenario. This result agrees with the minimax rate associated with
estimating a $(d+1)$-dimensional anisotropic density with respect to
the Hellinger distance (cf. \cite{Birge1986}). Interestingly, the
posterior automatically adapts to different vector of smoothness
levels. Noticeably, as in the isotropic case, the ambient dimension $p$
does not affect the rate except when it grows exponentially in $n$. It
will be interesting to generalize the result to allow anisotropic H\"{o}lder classes with noninteger smoothness levels as
in \cite{Shen2013}. Since approximation properties of tensor product of
B-splines are presently known only for integer smoothness in the
anisotropic case, we restrict to smoothness parameters to integer
values only.

\subsection{Deterministic predictor variables}\label{sec33}

Our method also applies for the case of deterministic predictors. In
order to obtain the posterior convergence rate, we first define the
empirical measure $\mathbb{P}_n^{\mathbf{X}} = n^{-1} \sum_{i=1}^n \delta
_{\mathbf{X}_i}$, in which $\delta_{\mathbf{X}_i}$ is a~point mass probability
measure at $\mathbf{X}_i$, based on the observations $\mathbf{X}_1,\ldots,\mathbf
{X}_n$. Then we define an empirical Hellinger distance on the space of
conditional densities by
%
\begin{equation}
\rho_{n}^2(f_1,f_2) = \iint
\bigl\{f_1^{1/2}(y|x_1,\ldots,x_p) -
f_2^{1/2}(y|x_1,\ldots,x_p) \bigr
\}^2 \,\mathrm{d}y \mathbb{P}_n^{\mathbf{X}}(\mathrm{d}x_1,
\ldots,\mathrm{d}x_p),
\end{equation}
where $f_1$ and $f_2$ are generic conditional densities of $Y$ on
$(0,1)$ given $\mathbf{X}$ in $(0,1)^p$. We can obtain the same posterior
contraction rates as the case of random predictors for both isotropic
and anisotropic cases.

\begin{theorem}\label{thm3}
Suppose that we have deterministic predictors $\mathbf{X}_1,\ldots,\mathbf{X}_n$ observed on $(0,1)^{p}$. Assume that the prior on the
conditional densities satisfies conditions \textup{(A1)}--\textup{(A4)}. If the true
conditional density satisfies \textup{(B1)}--\textup{(B3)}, then the posterior
distribution of $f$ contracts at $f_0$ at the rate $\varepsilon_n$ given
by \eqref{isotropicrate} with respect to $\rho_{n}$. If the true
conditional density satisfies \textup{(B2)}--\textup{(B4)}, and \textup{(A3)} is replaced by
\textup{(A3$'$)} in the prior specification, then the posterior distribution of
$f$ contracts at $f_0$ at the rate $\varepsilon_n$ given by \eqref{anisotropicrate} with respect to $\rho_{n}$.
\end{theorem}

\section{Numerical results}\label{secnum}

\subsection{Computation}
\label{subseccomp}
First, we ignore that we have a variable selection issue and pretend
that we know which $d$ predictors are relevant, say $\{m_1,\ldots,m_d\}
$. Thus, we may pretend that we are in the fixed dimensional setting
$p=d$ and all predictors are relevant. Then given the observations $(\mathbf
{X_i},Y_i)=(X_{i1},\ldots,X_{id};Y_i)$, $i=1,\ldots,n$,
the conditional likelihood
%
\begin{equation}
\label{new1}
L(\bolds{\eta},\mathbf{J}| \mathbf{X},Y) = \prod
_{i=1}^n \Biggl\{\sum_{j_0=1}^{J_0}
\sum_{j_{m_1}=1}^{J_{m_1}}\cdots\sum
_{j_{m_d}=1}^{J_{m_d}} \eta_{j_0,j_{m_1},\ldots,j_{m_d}} \bar
{B}_{j_0}(Y_i) \prod_{k=1}^{d}
B_{j_{m_k}}(X_{i, j_{m_k}}) \Biggr\}
\end{equation}
expands to $\sum_{\mathbf{s_1} \in\Lambda_{\mathbf{J}}} \cdots\sum_{\mathbf{s_n}
\in\Lambda_{\mathbf{J}}} \prod_{i=1}^n \eta_{\mathbf{s}_i} B_{\mathbf{s}_i}^{*}
(Y_i|\mathbf{X_i})$,
where $\Lambda_{\mathbf{J}} = \{1,\ldots,J_0\}\times\cdots\times\{
1,\ldots,\break J_{m_d}\}$ and $B_{\mathbf{s}}^{*}(y|\mathbf{x})$ is defined as $\bar
{B}_{s_0}(y) \prod_{k=1}^d B_{s_{m_k}}(x_{m_k})$ for every
$(d+1)$-dimensional index $\mathbf{s}$ taking value in $\Lambda_{\mathbf{J}}$
and $\mathbf{J}=(J_0,J_{m_1},\ldots,J_{m_d})\in\mathbb{N}^{d+1}$. Since
\[
\P(\mathbf{J},\bolds{\eta}|\mathbf{X},Y) \propto\P(\mathbf{X},Y|\mathbf{J},
\bolds{\eta}) \Pi _4(\bolds{\eta}|\mathbf{J})\tilde{
\Pi}_3(J_0)\prod_{k=1}^d
\tilde{\Pi}_3(J_{m_k}),
\]
the posterior mean of $f(y|\mathbf{x})$ at point is given by
%
\begin{equation}
\label{eqe111}
\frac{ \sum_{J_0=1}^{\infty} \cdots\sum_{J_{m_d}=1}^{\infty
} \tilde{\Pi}_3(J_0)\prod_{k=1}^d \tilde{\Pi}_3(J_{m_k}) \int f(y|\mathbf
{x},\bolds{\eta},\mathbf{J}) \Pi_4(\bolds{\eta}) L( \bolds{\eta},\mathbf{J}|\mathbf{X},Y )
\,\mathrm{d}\bolds{\eta} }{
 \sum_{J_0=1}^{\infty} \cdots\sum_{J_{m_d}=1}^{\infty}
\tilde{\Pi}_3(J_0)\prod_{k=1}^d \tilde{\Pi}_3(J_{m_k}) \int\Pi_4(\bolds
{\eta}) L( \bolds{\eta},\mathbf{J}|\mathbf{X},Y )\, \mathrm{d}\bolds{\eta}}.
\end{equation}

In view of \eqref{new1} and the form of $f(y|\mathbf{x})$ given by \eqref{eqe1000},
both numerator and denominator of \eqref{eqe111} involve
computing integrals of the form $I(\mathbf{s_1},\ldots,\mathbf{s_n}) = \int_{\bolds
{\eta} } \Pi(\bolds{\eta}) \prod_{k=1}^n \eta_{\mathbf{s_k}} \,\mathrm{d}\bolds{\eta}$. We
collect coefficients $\bolds{\eta}$ with the same index together to form
their powers and observe that, by (A4), coefficients whose index differ
anywhere except in the zeroth co-ordinate are independent, and the
collection of coefficients with the same last $d$ components are
Dirichlet distributed. In view of the conjugacy of the Dirichlet
functional form with respect to a multinomial function, $I(\mathbf{s}_1,\ldots,\mathbf{s}_n)$ can be written down in terms of products of
certain Dirichlet forms, and hence can be computed for any given $(\mathbf
{s}_1,\ldots,\mathbf{s}_n)$.
Therefore \eqref{eqe111} simplifies to
%
\begin{eqnarray}
&& \hspace*{-10pt}\frac{ \sum_{J_0=1}^{\infty} \cdots\sum_{J_d=1}^{\infty}
\tilde{\Pi}_3(J_0)\prod_{k=1}^d \tilde{\Pi}_3(J_{m_k}) \sum_{\mathbf{s_0}
\in \Lambda_{\mathbf{J}}} \cdots\sum_{\mathbf{s_n} \in \Lambda_{\mathbf{J}}}
I(\mathbf{s}_0,\ldots,\mathbf{s}_n) \prod_{i=0}^n B_{\mathbf{s}_i}^{*} (Y_i|\mathbf
{X_i}) }{
  \sum_{J_0=1}^{\infty} \cdots\sum_{J_d=1}^{\infty}
\tilde{\Pi}_3(J_0)\prod_{k=1}^d \tilde{\Pi}_3(J_{m_k}) \sum_{\mathbf{s}_1
\in \Lambda_{\mathbf{J}}} \cdots\sum_{\mathbf{s}_n \in \Lambda_{\mathbf{J}}}
I(\mathbf{s}_1,\ldots,\mathbf{s}_n) \prod_{i=1}^n B_{\mathbf{s}_i}^{*} (Y_i|\mathbf
{X}_i) },\nn
\\[-5pt]
\end{eqnarray}
where $(\mathbf{X_0},Y_0)$ stands for $(\mathbf{x},y)$.

Now, we take the variable selection into consideration. Suppose that
the proposed model size is $r$, which follows the prior distribution
$\Pi_1$. Given $r$, let the covariates $X_{m_1},\ldots,X_{m_r}$ enter
the model with probability $\Pi_2(m_1,\ldots,m_r|r)$. Define
\begin{eqnarray*}
W^0 (m_1,\ldots,m_r|\mathbf{X},Y) &= & \sum
_{J_0=1}^{\infty} \sum
_{J_{m_1}=1}^{\infty}\cdots\sum_{J_{m_r}=1}^{\infty}
\tilde{\Pi }_3(J_0) \prod_{k=1}^r
\tilde{\Pi}_3(J_{m_k})
\\
&&{}\times\sum_{\mathbf{s}_0 \in \Lambda_{\mathbf{J}}} \cdots\sum
_{\mathbf{s}_n
\in \Lambda_{\mathbf{J}}} I(\mathbf{s}_0,\ldots,\mathbf{s}_n)
\prod_{i=0}^n B_{\mathbf{s}_i}^{*}
(Y_i|\mathbf {X}_i),
\\
W^1 (m_1,\ldots,m_r|\mathbf{X},Y)&= &\sum
_{J_0=1}^{\infty} \sum
_{J_{m_1}=1}^{\infty}\cdots\sum_{J_{m_r}=1}^{\infty}
\tilde{\Pi }_3(J_0) \prod_{k=1}^r
\tilde{\Pi}_3(J_{m_k})
\\
&&{}\times\sum_{\mathbf{s}_1 \in \Lambda_{\mathbf{J}}} \cdots\sum
_{\mathbf{s}_n
\in \Lambda_{\mathbf{J}}} I(\mathbf{s}_1,\ldots,\mathbf{s}_n)
\prod_{i=1}^n B_{\mathbf{s}_i}^{*}
(Y_i,\mathbf {X}_i).
\end{eqnarray*}
Then the posterior mean of $f(y|\mathbf{x})$ is given by
%
\begin{equation}
\label{eqndraw}
\frac{ \sum_{r=1}^{\bar{r}} \Pi_1(r) \sum_{1 \leq m_1 <
\cdots< m_r \leq p} \Pi_2(m_1,\ldots,m_r|r) W^0(m_1,\ldots,m_r|\mathbf
{X},Y) }{
 \sum_{r=1}^{\bar{r}} \Pi_1(r) \sum_{1 \leq m_1 < \cdots<
m_r \leq p} \Pi_2(m_1,\ldots,m_r|r) W^1(m_1,\ldots,m_r|\mathbf{X},Y)}.
\end{equation}
Similar expressions can be obtained for other posterior moments, in
particular, for the posterior second moment and hence the posterior
variance. This estimate can be viewed as a kernel mixture estimator
whose kernel is determined jointly by selected covariates and
associated tensor product B-splines. Since a B-spline basis function
takes nonzero values only at $q$ intervals, the calculation of $W^0$
for a given $r$ involves $(J_{\max} - J_{\min} + 1)^{r+1}
q^{(r+1)(n+1)}$ terms if we restrict $J_0$ and each $J_{m_k}$,
$k=1,\ldots,r$, to take values between $J_{\min}$ and $J_{\max}$. Then
there will be $\sum_{r=1}^{\bar{r}} {p \choose r} (J_{\max} - J_{\min}
+ 1)^{r+1} q^{(r+1)(n+1)}$ terms in total. Instead of evaluating all
terms, we randomly sample a number of terms in both numerator and
denominator and take the associated average values. If we choose $q=1$,
then the prior can be viewed as a multivariate random histogram and the
number of terms in the expression for the posterior mean will reduce to
$\sum_{r=1}^{\bar{r}} {p \choose r} (J_{\max} - J_{\min} + 1)^{r+1}$,
although the resulting density estimate will be less smooth and the
rate adaptation property of the posterior distribution will apply only
to smoothness up to order 1. We shall make this choice in our
simulation studies
to save on computational cost in exchange of sacrificing some smoothness.



\subsection{Simulation results}\label{secsim}
In the following, we provide more details in the prior construction of
our model.
\begin{enumerate}[(C4)]
\item[(C0)] We choose $q = 1$, which leads to histogram basis functions
(Haar basis).
\item[(C1)] We assign a uniform prior on the model size ranging from
$2$ to $\bar{r} = 7$. 
%
\item[(C2)] The prior probability of $\gamma_k$ follows a Bernoulli
distribution with parameter $w_k$ for $0 \leq w_k \leq1$ and
$k=1,\ldots,p$.
The values of $w_k$ can depend on the marginal correlation between
$X_k$ and $Y$.
%
\item[(C3)] Given the model size $r$ chosen, we generate a
zero-truncated Poisson random variable~$K$ with mean $\lambda= 100$ and
then assign the integer part of $K^{1/(r+1)}$ to the number of
expansion terms $J$. We restrict $J$ between $4$ and $8$, that is,
$J_{\min} = 4$ and $J_{\max} = 8$. Then (A3) holds for $\kappa' = \kappa
'' = 1$. 
%
\item[(C4)] Given $J_0$, we let the vector $(\theta_{j_0,j_{m_1},\ldots
, j_{m_r}}\dvt j_0=1,
\ldots,J_0)$ have the uniform distribution over the $J_0$-simplex for
every feasible value of $\mathbf{j}$. Then condition (A4) is satisfied for
$a = 1$.
\end{enumerate}

We apply the MCMC-free calculation method described in Section~\ref{subseccomp} on the following two examples,
%
\begin{eqnarray}
\label{denreg1} Y|\mathbf{X} &\sim& \operatorname{Beta}\bigl(4X_1 +
3X_2^2, 10X_2\bigr),
\\
\label{denreg2} Y|\mathbf{X} &\sim& \operatorname{Beta}\bigl(5X_2
\exp(2X_1), 5X_3^2 + 3X_4
\bigr).
\end{eqnarray}
For each example, we generate $p$ covariates $X_1,\ldots,X_p$ uniformly
from $[0.05,0.95]$. In the computation of \eqref{eqndraw}, we randomly
draw $N^* = 100$ or $500$ terms in the sums of $W^1$ and $W^0$ for
every fixed choice of $m_1,\ldots,m_r$ and $r$. We compare our method
(rsp) with least-squares kernel conditional density estimation (L-S)
developed by \cite{Sugiyama2010}, where they use $L_1$-regularization
to select variables. Prediction errors under the $L_2$-loss and their
maximum standard errors associated with~$10$ Monte Carlo replications
are summarized in Tables~\ref{tdenreg1} and \ref{tdenreg2}.

\begin{table}[b]
\tablewidth=\textwidth
\caption{Simulation example 1: true density generated by
\protect\eqref{denreg1}}\label{tdenreg1}
\begin{tabular*}{\tablewidth}{@{\extracolsep{\fill}}lllllll@{}}
\hline
& \multicolumn{3}{l}{$n=100$}&\multicolumn{3}{l@{}}{$n=500$}\\[-4pt]
& \multicolumn{3}{l}{\hrulefill} & \multicolumn{3}{l@{}}{\hrulefill}\\
& L-S & rsp($N^*=100$) & rsp($N^*=500$) & L-S & rsp($N^*=100$) &
rsp($N^*=500$) \\ \hline
$p=5$ & $0.88$ & $0.76$ & $0.86$ & $0.81$ & $1.01$ & $1.01$ \\
$p=10$ & $1.04$ & $0.78$ & $0.81$ & $1.10$ & $1.10$ & $1.12$ \\
$p=50$ & $0.69$ & $0.66$ & $0.68$ & $0.73$ & $1.05$ & $1.04$ \\
$p=100$ & $0.69$ & $0.72$ & $0.70$ & $0.63$ & $0.76$ & $0.96$ \\
$p=500$ & $0.96$ & $0.78$ & $0.67$ & $0.86$ & $0.76$ & $0.80$ \\
$p=1000$ & $1.15$ & $0.59$ & $0.63$ & $1.26$ & $0.77$ & $0.95$ \\
max s.e. & $0.07$ & $0.06$ & $0.06$ & $0.08$ & $0.10$ & $0.10$ \\
\hline
\end{tabular*}
\end{table}

Compared with the least-squares method, our approach has a better
performance in most cases. Since we are directly sampling a fixed
number of terms from the sums in \eqref{eqndraw}, our prediction error
does not change too much with $p$, which makes the proposed method
outperform L-S when $p$ is large. Comparing the prediction errors under
the choices of $N^* = 100$ and $500$, we find that their performances
are quite close to each other. Hence direct sampling does not introduce
too much variability. We also carry out a sensitivity analysis by using
different parameter values in the prior distribution, for example, $\bar
{r}=6$, $\lambda= 50$ and $J_{\min} = 5$, $J_{\max}=10$. Similar
results are obtained. In practice, one may choose $\bar{r}$ as a
constant multiple (e.g., twice) of the possible maximal model size to
let all important covariates be included in the considered model with a
high probability. The Poisson mean parameter $\lambda$ in (C3) shall be
modified according to the choice of $\bar{r}$ to ensure that $\lambda
^{1/r}$ falls into an appropriate range, say, between $4$ and $20$.

\begin{table}
\caption{Simulation example 2: true density generated by
\protect\eqref{denreg2}} \label{tdenreg2}
\begin{tabular*}{\tablewidth}{@{\extracolsep{\fill}}lllllll@{}}
\hline
& \multicolumn{3}{l}{$n=100$}&\multicolumn{3}{l@{}}{$n=500$}\\[-4pt]
& \multicolumn{3}{l}{\hrulefill} &     \multicolumn{3}{l@{}}{\hrulefill}\\
& L-S & rsp($N^*=100$) & rsp($N^*=500$) & L-S & rsp($N^*=100$) &
rsp($N^*=500$) \\
\hline
$p=5$ & $0.76$ & $0.61$ & $0.64$ & $0.86$ & $0.74$ & $0.71$ \\
$p=10$ & $0.97$ & $0.66$ & $0.61$ & $1.00$ & $0.70$ & $0.71$ \\
$p=50$ & $0.69$ & $0.64$ & $0.61$ & $0.72$ & $0.71$ & $0.74$ \\
$p=100$ & $0.72$ & $0.67$ & $0.64$ & $0.74$ & $0.73$ & $0.72$ \\
$p=500$ & $0.95$ & $0.61$ & $0.71$ & $1.00$ & $0.73$ & $0.78$ \\
$p=1000$ & $1.26$ & $0.68$ & $0.66$ & $1.25$ & $0.69$ & $0.73$ \\
max s.e. & $0.08$ & $0.04$ & $0.05$ & $0.08$ & $0.06$ & $0.06$ \\
\hline
\end{tabular*}
\end{table}

\section{Proofs}\label{sec5}
\begin{pf*}{Proof of Theorem~\ref{thm}}
Note that the conditional density $f(y|\mathbf{x})$ is the same as the
joint density of $(\mathbf{X},Y)$ at $(\mathbf{x},y)$ with respect to the
dominating measure $\mu$ equal to the product of $G$ and the Lebesgue
measure. Further, the distance $\rho$ on the space of conditional
densities is equivalent to the Hellinger distance on the space of joint
densities with respect to $\mu$. Hence, in order to derive contraction
rate of the posterior distribution of the conditional density at a true
density $f_0(\cdot|\cdot)$, we need only to apply the standard result
on posterior convergence rate for (joint) densities given by Theorem~1
of \cite{Ghosal2000}. The required conditions characterizing the
posterior contraction rate $\varepsilon_n\to0$ can therefore be rewritten
in the present context as follows: there exists a sequence of subsets
$\mathcal{F}_n$ of the space of conditional densities, called sieves,
such that
%
\begin{eqnarray}\label{eqa112}
\Pi\bigl(\mathcal{F}_{n}^c\bigr) &\lesssim &  \exp\bigl
\{-8n\varepsilon_n^2\bigr\},
\\
\label{eqa111} \log D(\varepsilon_n,\mathcal{F}_{n}, \rho)  &\lesssim &  n
\varepsilon_n^2,
\\
\label{eqa1131} \Pi\bigl(f\dvt K(f_0,f) \leq\varepsilon_n^2,
V(f_0,f) \leq\varepsilon_n^2 \bigr) &\gtrsim  & \exp
\bigl\{- n \varepsilon_n^2\bigr\},
\end{eqnarray}
where $K(f_0,f) = \iint f_0(y|\mathbf{x}) \log(f_0(y|\mathbf{x})/f(y|\mathbf
{x})) \,\mathrm{d}y \, \mathrm{d}G(\mathbf{x})$ is the Kullback--Leibler divergence and $V(f_0,f)
= \iint f_0(y|\mathbf{x}) \log^2 (f_0(y|\mathbf{x})/f(y|\mathbf{x}))\, \mathrm{d}y  \,\mathrm{d}G(\mathbf
{x})$ is the Kullback--Leibler variation.
We define a sieve in the following way:
%
\begin{eqnarray}
\mathcal{F}_n &=&  \bigl\{h(
\mathbf{x},y|r;m_1,\ldots ,m_r;J_0,J_{m_1},
\ldots,J_{m_r};\bolds{\theta})\dvt r \leq\bar{r}_n,
\nonumber
\\[-8pt]
\label{eqsieve}
\\[-8pt]
\nonumber
&& \hspace*{4pt} 1 \leq m_1 < \cdots< m_r \leq p;
J_0,J_{m_1},\ldots,J_{m_r} \leq
\tilde{J}_n, \bolds{\theta} \in(\Delta_{J_0})^{\prod_{k=1}^r J_{m_k}}
\bigr\},
\end{eqnarray}
where $\tilde{J}_n = \lfloor(L J_n^*)^{(d+1)/(r+1)} (\log n)^{\kappa
/(r+1)} \rfloor$, $J_n^*$ and $\bar{r}_n$ are two sequences of number
going to infinity, $L$ and $\kappa$ are some fixed positive constants.
We shall choose the values of these numbers later.

We first verify \eqref{eqa112}. Note that $\Pi(\mathcal{F}_{n}^c) $ is
bounded by
\begin{eqnarray*}
&&  \Pi_1(r > \bar{r}_n) + \sum
_{r=1}^{\bar{r}_n} \sum_{1 \leq
m_1 < \cdots< m_r \leq p}
\Pi_3(J_{m_k} > \tilde{J}_n \mbox{ for some }k=1,\ldots,r |r,m_1,\ldots,m_r)
\\
&&\quad  \leq \exp\bigl\{-\exp\bigl(c_0 \bar{r}_n^{t_0}
\bigr) \bigr\} + \sum_{r=1}^{\bar{r}_n} {p \choose
r} r \tilde{\Pi}_3 ( J > \tilde{J}_n)
\\
&& \quad \leq \exp\bigl\{-\exp\bigl(c_0 \bar{r}_n^{t_0}
\bigr) \bigr\} + \bar{r}_n p^{\bar{r}_n} \sum
_{r=1}^{\bar{r}_n} \exp\bigl\{-c_2''
\tilde{J}_n^{r+1}(\log\tilde {J}_n)^{\kappa''}
\bigr\}
\\
&&\quad  \leq \exp\bigl\{-\exp\bigl(c_0 \bar{r}_n^{t_0}
\bigr) \bigr\} + \bar{r}_n^2 p^{\bar
{r}_n} \exp\bigl
\{-c_3 L^{d+1} \bigl(J_n^*\bigr)^{d+1}
(\log\tilde{J}_n)^{\kappa'' +
\kappa} \bigr\}
\\
&&\quad  =  \exp\bigl\{-\exp\bigl(c_0 \bar{r}_n^{t_0}
\bigr) \bigr\} + \exp\bigl\{2\log\bar{r}_n + \bar {r}_n
\log p - c_3 L^{d+1} \bigl(J_n^*
\bigr)^{d+1} (\log\tilde{J}_n)^{\kappa'' +
\kappa} \bigr\}
\\
&&\quad  \leq  \exp\bigl(-b n \varepsilon_n^2\bigr)
\end{eqnarray*}
for any $b > 0$ and some constant $c_3 > 0$ provided $L$ is chosen
sufficiently large and the following relations hold
%
\begin{eqnarray}
&& \log\bar{r}_n \lesssim n \varepsilon_n^2,
\qquad \bar{r}_n \log p \lesssim n\varepsilon_n^2,
\nonumber
\\[-8pt]
\label{rel1}\\[-8pt]
\nonumber
&& \bigl(J_n^*\bigr)^{d+1} (\log n )^{\kappa+ \kappa''} \gtrsim n
\varepsilon_n^2, \qquad \exp\bigl(c_0
\bar{r}_n^{t_0}\bigr) \gtrsim n\varepsilon_n^2.
\end{eqnarray}

Now we bound the covering number $D(\varepsilon_n, \mathcal{F}_n, \rho)$
using the relation $D(\varepsilon_n, \mathcal{F}_n, \rho) \leq D(\varepsilon
_n^2, \mathcal{F}_n, \|\cdot\|_1)$, where $\|\cdot\|_1$ stand for the
$L_1$-distance on the space of conditional densities given by
\[
\|f_1-f_2\|_1=\iint \bigl\llvert
f_1(y|\mathbf{x})-f_2(y|\mathbf{x})\bigr\rrvert \,\mathrm{d}y \,\mathrm{d}G(
\mathbf{x})\le\sup_{\mathbf{x}} \int \bigl|f_1(y|
\mathbf{x})-f_2(y|\mathbf{x})\bigr|\,\mathrm{d}y.
\]
We split $\mathcal{F}_n$ in layers corresponding to different $r$,
different $m_1,\ldots,m_r$ and different $J_0,J_{m_1},\allowbreak \ldots,J_{m_r}$:
\[
\mathcal{F}_n = \bigcup_{r=1}^{\bar{r}_n}
\bigcup_{1 \leq m_1 < \cdots
< m_r \leq p,} \bigcup_{1 \leq J_0,J_{m_1},\ldots,J_{m_r} \leq J_n^*}
\mathcal{F}_{n;r;m_1,\ldots,m_r;J_0,J_{m_1},\ldots,J_{m_r}}.
\]
For any given $r$, $m_1,\ldots,m_r$, $J_0,J_{m_1},\ldots,J_{m_r}$,
consider $\bolds{\theta},\bolds{\theta}' \in(\Delta_{J_0})^{\prod_{k=1}^d
J_{m_k}}$. We can write
$\bolds{\theta} = (\bolds{\theta}_{j_{m_1},\ldots,j_{m_r}}\dvt 1 \leq
j_{m_1},\ldots,j_{m_r} \leq J_n^*)$, $\bolds{\theta}' = (\bolds{\theta
}'_{j_{m_1},\ldots,j_{m_r}}\dvt 1 \leq j_{m_1},\ldots,j_{m_r} \leq J_n^*)$
where $\bolds{\theta}_{j_{m_1},\ldots,j_{m_r}} = (\theta
_{j_0,j_{m_1},\ldots,j_{m_r}}\dvt 1 \leq j_0 \leq J_n^*)$ and $\bolds{\theta
}_{j_{m_1},\ldots,j_{m_r}}' = (\theta_{j_0,j_{m_1},\ldots,j_{m_r}}'\dvt 1
\leq j_0 \leq J_n^*)$. Then
\begin{eqnarray*}
&& \bigl\llVert h(\mathbf{x},y|r; m_1,\ldots,m_r;J_0,J_{m_1},
\ldots ,J_{m_r};\bolds{\theta}) - h\bigl(\mathbf{x},y|r;
m_1,\ldots,m_r;J_0,J_{m_1},\ldots
,J_{m_r};\bolds{\theta}'\bigr) \bigr\rrVert
_1
\\[-1pt]
&&\quad  \leq \sup_{\mathbf{x}} \sum_{j_0=1}^{J_0}
\sum_{j_{m_1}=1}^{J_{m_1}} \cdots\sum
_{j_{m_r}=1}^{J_{m_r}} \bigl| \theta_{j_0,j_{m_1},\ldots,j_{m_r}} -
\theta_{j_0,j_{m_1},\ldots,j_{m_r}}' \bigr| B_{j_{m_1}}(x_{m_1})
\cdots B_{j_{m_r}}(x_{m_r})
\\[-1pt]
&&\quad  \leq \max_{j_{m_1},\ldots,j_{m_r}} \bigl\|\bolds{\theta}_{j_{m_1},\ldots
,j_{m_r}} - \bolds{
\theta}_{j_{m_1},\ldots,j_{m_r}}' \bigr\|_1
\end{eqnarray*}
since the collection $B_j(x)$s add up to 1 for any $x$. Using the fact
that $D(\varepsilon,\Delta_d,\|\cdot\|_1) \leq(3/\varepsilon)^d$, we obtain
%
\begin{eqnarray}
D\bigl(\varepsilon_n^2,\mathcal{F}_{n;r;m_1,\ldots,m_r;J_0,J_{m_1},\ldots
,J_{m_r}},\|\cdot
\|_1\bigr)
& \leq  &  \prod_{1 \leq J_{m_1},\ldots,J_{m_r} \leq\tilde{J}_n} D\bigl(
\varepsilon_n^2,\Delta_{J_0},\|\cdot
\|_1\bigr) \nn
\\[-1pt]
& \leq &  \prod_{1 \leq J_{m_1},\ldots,J_{m_r} \leq\tilde{J}_n} \biggl(\frac{3}{\varepsilon_n^2}
\biggr)^{J_0}
\nonumber
\\[-9pt]
\\[-9pt]
\nonumber
&  = &  \biggl(\frac{3}{\varepsilon_n^2} \biggr)^{\tilde{J}_n^{r+1}}
\\[-1pt]
&  \leq &  \biggl(\frac{3}{\varepsilon_n^2} \biggr)^{L^{d+1} (J_n^*)^{d+1}
(\log n)^{\kappa}}.\nn
\end{eqnarray}
Therefore,\vspace*{-2pt}
\begin{eqnarray*}
 D(\varepsilon_n,\mathcal{F}_n,\rho) &\leq& D\bigl(
\varepsilon_n^2,\mathcal{F}_n,\|\cdot
\|_1\bigr) \nn
\\[-1.5pt]
& \leq& \sum_{r=1}^{\bar{r}_n} \sum
_{1 \leq m_1 < \cdots< m_r \leq
p,} \sum_{1 \leq J_0,J_{m_1},\ldots,J_{m_r} \leq\tilde{J}_n} \biggl(
\frac{3}{\varepsilon_n^2} \biggr)^{L^{d+1} (J_n^*)^{d+1} (\log n)^{\kappa
}} \nn
\\[-1.5pt]
& \leq& \sum_{r=1}^{\bar{r}_n} {p \choose r}
\tilde{J}_n^{r+1} \biggl(\frac{3}{\varepsilon_n^2}
\biggr)^{L^{d+1} (J_n^*)^{d+1} (\log n)^{\kappa
}} \nn
\\[-1.5pt]
& \le& \sum_{r=1}^{\bar{r}_n} p^r
L^{d+1} \bigl(J_n^*\bigr)^{d+1} (\log
n)^{\kappa} \biggl(\frac{3}{\varepsilon_n^2} \biggr)^{L^{d+1} (J_n^*)^{d+1}
(\log n)^{\kappa}} \nn
\\[-1.5pt]
& \lesssim& \bar{r}_n p^{\bar{r}_n} \bigl(J_n^*
\bigr)^{d+1} (\log n)^{\kappa} \exp \bigl\{L^{d+1}
\bigl(J_n^*\bigr)^{d+1} (\log n)^{\kappa} \log\bigl(3/
\varepsilon _n^2\bigr) \bigr\} \nn
\\[-1.5pt]
& = & \exp\bigl\{\log\bar{r}_n + \bar{r}_n \log p +
(d+1) \log J_n^* + \kappa\log(\log n) \nn
\\
&& \hspace*{18pt}{}+ L^{d+1} \bigl(J_n^*\bigr)^{d+1} (\log
n)^{\kappa} \log\bigl(3/\varepsilon _n^2\bigr)\bigr\}
\nn
\\[-1.5pt]
& \leq& \exp\bigl\{c_4 \bigl(J_n^*
\bigr)^{d+1} (\log n)^{\kappa+ 1} + \bar{r}_n \log p \bigr
\} \nn
\end{eqnarray*}
for some $c_4 > 0$. Thus it suffices to have the following relations
%
\begin{equation}
\label{rel2}
\bigl(J_n^*\bigr)^{d+1} (\log
n)^{\kappa+ 1} \lesssim n \varepsilon_n^2,
\qquad  \bar{r}_n \log p \lesssim n \varepsilon_n^2.
\end{equation}

For \eqref{eqa1131}, in order to lower bound the prior concentration
probability around $f_0(y|\mathbf{x})$, we shall restrict to the oracle
model consisting of $d$ true covariates $X_{m_1^0},\ldots,X_{m_d^0}$.
By Lemma~\ref{lemmal00}, there exists $\bolds{\theta_0} = (\theta
_{j_0,j_{m_1^0},\ldots,j_{m_d^0}}^0\dvt 1 \leq j_0, j_{m_1^0},\ldots
,j_{m_d^0} \leq J_n^*)$ such\vspace*{-2pt} that
%
\begin{equation}
\label{eqee0}
\sup_{\mathbf{x},y} \bigl\llvert f_0(y|
\mathbf{x} ) - h\bigl(\mathbf{x},y|d;m_1^0,\ldots
,m_d^0;J_n^*,\ldots,J_n^*;\bolds{
\theta_0} \bigr) \bigr\rrvert \lesssim \bigl(J_n^*
\bigr)^{-\beta} \leq\varepsilon_n.
\end{equation}
Now for every $(j_{m_1^0},\ldots,j_{m_d^0})$, we define $\bolds{\theta
}_{j_{m_1^0},\ldots,j_{m_d^0}}^0 = (\theta_{j_0,j_{m_1^0},\ldots
,j_{m_d^0}}^0\dvt 1 \leq j_0 \leq J_n^*) \in\Delta_{J_n^*}$. Then $\bolds
{\theta}_0$ can be written by $\bolds{\theta}_0= (\bolds{\theta
}_{j_{m_1^0},\ldots,j_{m_d^0}}^0\dvt 1 \leq j_{m_1^0},\ldots,j_{m_d^0}
\leq J_n^*)$. We consider $\bolds{\theta} = (\bolds{\theta}_{j_{m_1^0},\ldots
,j_{m_d^0}}\dvt 1 \leq j_{m_1^0},\ldots,j_{m_d^0} \leq J_n^*) $ and $\bolds
{\theta}_{j_{m_1^0},\ldots,j_{m_d^0}} =(\theta_{j_0,j_{m_1^0},\ldots
,j_{m_d^0}}\dvt 1 \leq j_0 \leq J_n^*) \in\Delta_{J_n^*} $.\vspace*{-2pt} If
%
\begin{equation}
\max_{1 \leq j_{m_1^0},\ldots,j_{m_d^0} \leq J_n^*} \bigl\|\bolds{\theta }_{j_{m_1^0},\ldots,j_{m_d^0}} - \bolds{
\theta}_{j_{m_1^0},\ldots
,j_{m_d^0}}^0 \bigr\|_1 \leq\varepsilon,
\end{equation}
then\vspace*{-2pt}
\begin{eqnarray*}
&& \bigl| h\bigl(\mathbf{x},y|d;m_1^0,
\ldots,m_d^0;J_n^*,\ldots,J_n^*;
\bolds {\theta} \bigr) - h\bigl(\mathbf{x},y|d;m_1^0,
\ldots,m_d^0;J_n^*,\ldots,J_n^*;
\bolds {\theta}_0 \bigr) \bigr|
\\
&& \quad  \leq  \sum_{j_0=1}^{J_n^*} \sum
_{j_{m_1^0}=1}^{J_n^*} \cdots\sum_{j_{m_d^0}=1}^{J_n^*}
\bigl| \theta_{j_0,j_{m_1^0},\ldots,j_{m_d^0}} - \theta_{j_0,j_{m_1^0},\ldots,j_{m_d^0}}^0 \bigr|
\bar{B}_{j_0}(y) B_{j_{m_1^0}} (x_{m_1^0}) \cdots
B_{j_{m_d^0}} (x_{m_d^0}).
\end{eqnarray*}
Since $0 \leq B_j(x) \leq1$ and $ 0 \leq\bar{B}_j(y) \leq J_n^*$ for
any $j$, we have
%
\begin{eqnarray}
&& \sup_{\mathbf{x},y} \bigl| h\bigl(\mathbf{x},y|d;m_1^0,
\ldots,m_d^0;J_n^*,\ldots ,J_n^*,
\bolds{\theta} \bigr) - h\bigl(\mathbf{x},y|d;m_1^0,
\ldots,m_d^0;J_n^*,\ldots ,J_n^*,
\bolds{\theta}_0 \bigr) \bigr|
\nonumber
\\[-8pt]
\label{eqee1}
\\[-8pt]
\nonumber
&& \quad \leq\bigl(J_n^*\bigr)^{d+1} \max_{1 \leq j_{m_1^0},\ldots,j_{m_d^0} \leq J_n^*}
\bigl\|\bolds{\theta}_{j_{m_1^0},\ldots,j_{m_d^0}} - \bolds{\theta }_{j_{m_1^0},\ldots,j_{m_d^0}}^0
\bigr\|_1 \leq\varepsilon_n
\end{eqnarray}
provided that
%
\begin{equation}
\label{tobound}
\bigl\|\bolds{\theta}_{j_{m_1^0},\ldots,j_{m_d^0}} - \bolds{\theta
}_{j_{m_1^0},\ldots,j_{m_d^0}}^0 \bigr\|_1 \leq\bigl(J_n^*
\bigr)^{-(d+1)} \varepsilon_n \qquad \mbox{for all }j_{m_1^0},
\ldots,j_{m_d^0}.
\end{equation}
To simplify the notation, we denote $h(\mathbf{x},y|d;m_1^0,\ldots
,m_d^0;J_n^*,\ldots,J_n^*;\bolds{\theta} )$ by $f_{\bolds{\theta}}$.
Combining \eqref{eqee0} and \eqref{eqee1}, we have the desired
approximation $\sup_{\mathbf{x},y} |f_0(y|\mathbf{x}) - f_{\bolds{\theta}}(y|\mathbf
{x}) | \leq2\varepsilon_n$.

Using condition (B3), $\inf f_{\bolds{\theta}} \geq\inf f_0 - \|f_0 -
f_{\bolds{\theta}}\|_{\infty} \geq\underline{m}_0/2$ given that ${\varepsilon
}_n$ is sufficiently small. This implies that $\|f_0/f_{\bolds{\theta}}\|
_\infty\le2\|f_0\|_\infty/\underline{m}_0<\infty$ since $f_0$ can be\vspace*{1pt}
regarded as a fixed continuous function on the compact set
$[0,1]^{d+1}$. Hence, for every $f_{\bolds{\theta}}$ satisfying $\|f_{\bolds
{\theta}} - f_0\|_{\infty} \leq2\varepsilon_n$,
%
\begin{equation}
\rho^2(f_0, f_{\bolds{\theta}}) = \int\frac{|f_0(y|\mathbf{x})-f_{\bolds{\theta
}}(y|\mathbf{x})|^2}{ (f_0^{1/2}(y|\mathbf{x})+f_{\bolds{\theta}}^{1/2}(y|\mathbf
{x}) )^{2}}\,\mathrm{d}y
\,\mathrm{d}G(\mathbf{x}) \leq\frac{1}{\underline{m}_0} \|f_0 - f_{\bolds{\theta}}
\|_{\infty}^2 \lesssim \varepsilon_n^2.
\end{equation}
Therefore, in view of Lemma 8 of \cite{Ghosal2007}, we obtain
%
\begin{eqnarray}
K(f_0,f_{\bolds{\theta}}) & \leq& 2
\rho^2(f_0,f_{\bolds{\theta}}) \biggl\| \frac{f_0}{f_{\bolds{\theta}}}
\biggr\|_{\infty} \lesssim\varepsilon_n^2 ,
\nonumber
\\[-8pt]
\label{eqe401}
\\[-8pt]
\nonumber
V(f_0,f_{\bolds{\theta}}) & \lesssim& \rho^2(f_0,f_{\bolds{\theta}})
\biggl( 1 + \biggl\| \frac{f_0}{f_{\bolds{\theta}}} \biggr\|_{\infty} \biggr)^2
\lesssim \varepsilon_n^2.
\end{eqnarray}
Thus, it suffices to lower bound the prior probability of the event in
\eqref{tobound}, which is
\begin{eqnarray*}
&& \Pi_1(d) \times\Pi_2 \bigl(\bigl
\{m_1^0,\ldots,m_d^0\bigr\} |r=d
\bigr) \times \bigl\{\tilde{\Pi}_3\bigl(J_n^*\bigr)
\bigr\}^{d+1}  \nn
\\
&&\qquad  {}\times\prod_{1 \leq j_{m_1^0},\ldots,j_{m_d^0} \leq J_n^*} \tilde {\Pi}_4
\bigl(\bigl\|\bolds{\theta}_{j_{m_1^0},\ldots,j_{m_d^0}} - \bolds{\theta }_{j_{m_1^0},\ldots,j_{m_d^0}}^0
\bigr\|_1 \leq\bigl(J_n^*\bigr)^{-(d+1)} \varepsilon
_n \bigr) \nn
\\
&& \quad \gtrsim  \frac{1}{{p \choose d}} \exp \bigl\{-(d+1) c_2'
\bigl(J_n^*\bigr)^{d+1} \bigl(\log J_n^*
\bigr)^{\kappa'} \bigr\} \times\exp \biggl\{ -\bigl(J_n^*
\bigr)^d c_5 J_n^* \log\frac{(J_n^*)^{d+1}}{\varepsilon_n}
\biggr\} \nn
\end{eqnarray*}
for some constant $c_5 > 0$ by the small ball probability estimates of
a Dirichlet distribution in Lemma~6.1 of \cite{Ghosal2000}. As long as
$J_n^*$ and $\varepsilon_n^{-1}$ are powers of $n$ within slowly varying
factors, the last expression can be bounded below by $\exp\{-d\log p -
c_6 (J_n^*)^{d+1} (\log n)^{\kappa'} \} $ for some $c_6 > 0$. Hence in
order to obtain \eqref{eqa1131}, it suffices to have the following
relationships:
%
\begin{equation}
\label{rel3}
\bigl(J_n^*\bigr)^{-\beta} \lesssim
\varepsilon_n, \qquad \log p \lesssim n \varepsilon_n^2, \qquad
\bigl(J_n^*\bigr)^{d+1} (\log n)^{\kappa' } \lesssim n
\varepsilon_n^2.
\end{equation}
We can determine the rate $\varepsilon_n$ as the smallest sequence of
numbers that satisfies \eqref{rel1}, \eqref{rel2} and~\eqref{rel3},
that is,
%
\begin{eqnarray}
\varepsilon_n  &=&  \max \bigl\{n^{-(1-\alpha)/2} (\log n)^{1/(2t_0)},
n^{-\beta/(2\beta+d+1)} (\log n)^{\kappa' \beta/(2\beta+ d+1)} \bigr\},
\nonumber
\\[-8pt]
\\[-8pt]
\nonumber
J_n^*  &=& \bigl\lfloor\bigl(n\varepsilon_n^2
\bigr)^{1/(d+1)} (\log n)^{1/(d+1)} \bigr\rfloor+ 1,
\end{eqnarray}
and $\kappa= \kappa' - \kappa''$, $\bar{r}_n = L' (\log n)^{1/t_0}$
for some sufficiently large $L'$ provided that the condition $\exp\{\exp
(c_0 \bar{r}_n^{t_0})\} \gg n\varepsilon_n^2$ is satisfied.
\end{pf*}

\begin{pf*}{Proof of Theorem~\ref{thm2}}
The proof essentially follows the outline given in Theorem~\ref{thm}
except for two main differences. First, we shall need to use different
$J_0,J_{m_1^0},\ldots,J_{m_d^0}$ due to the approximation result by
Lemma~\ref{lemmal2}. In particular, we need $J_0^{-\beta_0} \asymp
J_{m_1^0}^{-\beta_1} \asymp\cdots\asymp J_{m_d^0}^{-\beta_d} $. This
will sightly change our definition of the sieve and the calculation of
the prior concentration. The second difference is that we now have a
dependent prior distribution in (A3$'$), which will change the
calculation of the prior concentration rate.

We define a new sieve as follows
%
\begin{eqnarray}
\mathcal{F}_n &=&  \Biggl\{h(
\mathbf{x},y|r;m_1,\ldots ,m_r;J_0,J_{m_1},
\ldots,J_{m_r};\bolds{\theta})\dvt r \leq\bar{r}_n,
\nonumber
\\[-8pt]
\label{eqsieve1}
\\[-8pt]
\nonumber
&&\hspace*{6pt} 1 \leq m_1 < \cdots< m_r \leq p; J_0
\prod_{k=1}^r J_{m_k} \leq
\tilde{J}_n^{r+1}, \bolds{\theta} \in\Delta_{J_0}^{\prod_{k=1}^r J_{m_k}}
\Biggr\},
\end{eqnarray}
where $\tilde{J}_n = \lfloor(L J_n^*)^{(d+1)/(r+1)} (\log n)^{\kappa
/(r+1)} \rfloor$, $J_n^*$ and $\bar{r}_n$ are two sequences of number
going to infinity, and $L$ and $ \kappa$ are some fixed positive
constants. We shall choose the values of these numbers later.

We first verify \eqref{eqa112}. It follows that
%
\begin{equation}
\label{to11}
\Pi\bigl(\mathcal{F}_{n}^c\bigr) \leq
\Pi_1(r > \bar{r}_n) + \sum_{r=1}^{\bar
{r}_n}
\sum_{1 \leq m_1 < \cdots< m_r \leq p} \Pi_3 \Biggl(J_0
\prod_{k=1}^r J_{m_k} >
\tilde{J}_n^{d+1} \Biggr).
\end{equation}
Note that the joint distribution of $(J_0,J_{m_1},\ldots,J_{m_r})$
depends only on the value of their product $J_0 \prod_{k=1}^r J_{m_k}$.
Let $t$ be a given integer and let $N_t$ stand for the number of ways
one can choose $\{J_0,J_{m_1},\ldots,J_{m_r}\}$ such that $J_0 \prod_{k=1}^r J_{m_k} = t$. Then
\[
\Pi_3 \Biggl(J_0 \prod_{k=1}^r
J_{m_k} = t \Biggr) \leq N_t \exp \bigl\{
-c_2''' t (\log
t)^{\kappa''} \bigr\}
\]
for some $c_2''' > 0$. Clearly $N_t \leq t^{r+1}$. Thus,
\[
\Pi_3 \Biggl(J_0 \prod_{k=1}^r
J_{m_k} = t \Biggr) \leq \exp \bigl\{(r+1) \log t -
c_2''' t (\log
t)^{\kappa''} \bigr\} \leq\exp\bigl\{- c_7 t (\log
t)^{\kappa''}\bigr\}
\]
for some $c_7 > 0$ provided that $ t (\log t)^{\kappa''} \gg(r+1) \log
t$, which is satisfied if $ t \gg r$ since $ \kappa'' \geq1$. Note
that the distribution of the product $J_0 \prod_{k=1}^r J_{m_k}$ has a
better-than-geometric tail starting from a large multiple of $r$, and hence
\[
\Pi_3\Biggl(J_0 \prod_{k=1}^r
J_{m_k} \geq t\Biggr) \leq\exp(-c_8 t \log t)
\]
for some $c_8 > 0$. In the sieve, we choose the cut-off $ \tilde
{J}_n^{r+1}$ which is clearly of order higher than $r$ and hence the
requirement is met. As a result, the second term in \eqref{to11} is
bounded by
\begin{eqnarray*}
&& \sum_{r=1}^{\bar{r}_n} p^r
\exp\bigl\{-c_8 \tilde{J}_n^{r+1} \log
\tilde{J}_n^{r+1} \bigr\}
\\
&& \quad  \leq  \sum_{r=1}^{\bar{r}_n} \exp\bigl\{r \log
p - c_8' L^{d+1} \bigl(J_n^*
\bigr)^{d+1} (\log n)^{\kappa+ \kappa''} \bigr\}
\\
&& \quad  \leq \exp\bigl\{\log\bar{r}_n + \bar{r}_n \log p -
c_8' L^{d+1} \bigl(J_n^*
\bigr)^{d+1} (\log n)^{\kappa+ \kappa''} \bigr\}
\end{eqnarray*}
for some $c_8' > 0$, which is of the same form of the corresponding
bound for the isotropic case, and that $L$ can be chosen sufficiently
large. Thus, relation \eqref{rel1} is obtained.

The calculation of entropy proceeds in the same way as in the isotropic
case. We split $\mathcal{F}_n$ into layers following the same
definition. Then
\[
D\bigl(\varepsilon_n^2,\mathcal{F}_{n;r;m_1,\ldots,m_r;J_0,J_{m_1},\ldots
,J_{m_r}},\|\cdot
\|_1\bigr) \leq \biggl(\frac{3}{\varepsilon_n^2} \biggr)^{J_0 \prod_{k=1}^r J_{m_k}} \leq
\biggl(\frac{3}{\varepsilon_n^2} \biggr)^{\tilde{J}_n^{r+1}}
\]
and the remaining calculations are identical, which give entropy
estimates of the sieve as in the isotropic case and hence relation
\eqref{rel2} is obtained.

Now we estimate the prior concentration rate. Consider the oracle model
given by $(d;m_1^0,\ldots, m_d^0; J_{n,0}^*,J_{n,1}^*,\ldots
,J_{n,d}^*)$, where
\[
\bigl(J_{n,0}^*\bigr)^{-\beta_0} \asymp\bigl(J_{n,1}^*
\bigr)^{-\beta_1} \asymp\cdots \asymp \bigl(J_{n,d}^*
\bigr)^{-\beta_d} \leq\varepsilon_n.
\]

By Lemma~\ref{lemmal2}, there exists $\bolds{\theta}_0 = (\theta
_{j_0,j_{m_1^0},\ldots,j_{m_d^0}}\dvt 1 \leq j_0 \leq J_{n,0}^*,1\leq
j_{m_k^0} \leq J_{n,k}^*,k=1,\ldots,d)$ such that
\[
\sup_{\mathbf{x},y} \bigl\llvert f_0(y|\mathbf{x}) - h
\bigl(\mathbf{x},y|d;m_1^0,\ldots ,m_d^0;J_{n,0}^*,
\ldots,J_{n,d}^*;\bolds{\theta}_0\bigr) \bigr\rrvert
\lesssim\sum_{k=0}^d \bigl(J_{n,k}^*
\bigr)^{-\beta_k} \lesssim\varepsilon_n.
\]
Given $j_{m_1^0},\ldots,j_{m_d^0}$, define $\bolds{\theta
}_{j_{m_1^0},\ldots,j_{m_d^0}}^0 = (\theta_{j_0,j_{m_1^0},\ldots
,j_{m_d^0}}^0\dvt 1 \leq j_0 \leq J_{n,0}^*) \in\Delta_{J_{n,0}^*}$. Then
$\bolds{\theta}_0 \in(\Delta_{J_{n,0}^*})^{\prod_{k=1}^d J_{n,k}^*}$. Let
$\bolds{\theta} \in(\Delta_{J_{n,0}^*})^{\prod_{k=1}^d J_{n,k}^*}$ and be
represented by $(\bolds{\theta}_{j_{m_1^0},\ldots,j_{m_d^0}}\dvt 1 \leq
j_{m_k^0} \leq J_{n,k}^*, k=1,\ldots,d)$. Then as before,
\begin{eqnarray*}
&& \bigl| h\bigl(\mathbf{x},y|d;m_1^0,
\ldots,m_d^0;J_{n,0}^*,\ldots
,J_{n,d}^*,\bolds{\theta} \bigr) - h\bigl(\mathbf{x},y|d;m_1^0,
\ldots ,m_d^0;J_{n,0}^*,\ldots,J_{n,d}^*,
\bolds{\theta_0} \bigr) \bigr|
\\
&& \quad  \leq  \sum_{j_0=1}^{J_n^*} \sum
_{j_{m_1^0}=1}^{J_n^*} \cdots\sum_{j_{m_d^0}=1}^{J_n^*}
\bigl| \theta_{j_0,j_{m_1^0},\ldots,j_{m_d^0}} - \theta_{j_0,j_{m_1^0},\ldots,j_{m_d^0}}^0 \bigr|
\bar{B}_{j_0}(y) B_{j_{m_1^0}}^*(x_{m_1^0}) \cdots
B_{j_{m_d^0}}^*(x_{m_d^0})
\\
&& \quad  \leq J_{n,0}^* \prod_{k=1}^d
J_{n,k}^* \max_{\stackrel{1 \leq
j_{m_k^0} \leq J_{n,k}^*,}{k=1,\ldots,d}} \bigl\|\bolds{\theta
}_{j_{m_1^0},\ldots,j_{m_d^0}} - \bolds{\theta}_{j_{m_1^0},\ldots
,j_{m_d^0}}^0 \bigr\|_1
\\
&& \quad  \leq  \bigl(J_n^*\bigr)^{d+1} \max_{\stackrel{1 \leq j_{m_k^0} \leq
J_{n,k}^*,}{k=1,\ldots,d}}
\bigl\|\bolds{\theta}_{j_{m_1^0},\ldots,j_{m_d^0}} - \bolds{\theta}_{j_{m_1^0},\ldots,j_{m_d^0}}^0
\bigr\|_1,
\end{eqnarray*}
where $J_n^* = \lfloor(\prod_{k=0}^d J_{n,k}^*)^{1/(d+1)} \rfloor+ 1$
is the smallest integer greater than the geometric mean of
$J_{n,0}^*,\ldots,J_{n,d}^*$. Thus it suffices to lower bound
\begin{eqnarray*}
&& \Pi_1(d) \times\Pi_2 \bigl(\bigl
\{m_1^0,\ldots,m_d^0\bigr\} |r=d
\bigr) \times\Pi_3\bigl(J_{n,0}^*,\ldots,J_{n,d}^*
\bigr)
\\
&& \quad {}\times\prod_{1 \leq j_0 \leq J_{n,0}^*, 1 \leq j_{m_k^0} \leq
J_{n,k}^*, k=1,\ldots,d} \tilde{\Pi}_4
\bigl(\bigl\|\bolds{\theta }_{j_{m_1^0},\ldots,j_{m_d^0}} - \bolds{\theta}_{j_{m_1^0},\ldots
,j_{m_d^0}}^0
\bigr\| \leq\bigl(J_n^*\bigr)^{-(d+1)} \varepsilon_n
\bigr).
\end{eqnarray*}
Since the other factors are as before, it suffices to look at the third
factor only, whose lower bound is given by
\[
\exp \Biggl\{ - c_2'' \Biggl(\,\prod
_{k=0}^d J_{n,k}^* \Biggr)
\Biggl(\, \sum_{k=0}^d \log
J_{n,k}^* \Biggr)^{\kappa} \Biggr\} \gtrsim\exp \bigl\{ -
c_9 \bigl(J_n^*\bigr)^{d+1} \bigl(\log
J_n^*\bigr)^{\kappa} \bigr\}
\]
for some constant $c_9 > 0$, which is identical with the corresponding
expression for the isotropic case. Thus we need
%
\begin{eqnarray}
&& \bigl(J_{n,k}^*\bigr)^{-\beta_k} \lesssim
\varepsilon_n,\qquad  \log p \lesssim n \varepsilon _n^2,
\qquad \bigl(J_n^*\bigr)^{d+1} (\log n)^{\kappa'} \lesssim n
\varepsilon_n^2,
\nonumber
\\[-8pt]
\label{rel31}\\[-8pt]
\nonumber
&&  J_n^* \asymp \Biggl(\,\prod
_{k=0}^d J_{n,k}^*
\Biggr)^{1/(d+1)}.
\end{eqnarray}
Combining\vspace*{1pt} \eqref{rel1}, \eqref{rel2} and \eqref{rel31}, we can choose
$\kappa= \kappa' - \kappa''$, $\bar{r}_n$ as a large multiple of
$(\log n)^{1/t_0}$, $J_{n,k}^* = \varepsilon_n^{-1/\beta_k}$ and
\[
\varepsilon_n = \max \bigl\{n^{-\beta^*/(2\beta^* + d + 1)} (\log n)^{\kappa
' \beta^*/(2\beta^* + d + 1)},
n^{-(1-\alpha)/2} (\log n)^{1/(2t_0)} \bigr\},
\]
where $\beta^* = (d+1)  (\sum_{k=0}^{d} \beta_k^{-1} )^{-1}$
is the harmonic mean of $(\beta_0,\beta_1,\ldots,\beta_d)$.
\end{pf*}

\begin{pf*}{Proof of Theorem~\ref{thm3}}
Note that the distance $\rho_{n}$ on the space of conditional densities
mathematically can be expressed as the Hellinger distance on the space
of joint densities with respect to the dominating measure $\mu$, which
is the product of $\mathbb{P}_n^{\mathbf{X}}$ and the Lebesgue measure.
This is notwithstanding the fact that the predictor variables are
actually deterministic. We only need to replace $G(\cdot)$ by $\mathbb
{P}_n^{\mathbf{X}}$ in the definitions of Kullback--Leibler divergence and
Kullback--Leibler variation in~\eqref{eqa1131}. The rest of arguments
proceed exactly in the same way as in Theorems \ref{thm} and \ref{thm2}
for the isotropic and the anisotropic cases respectively.
\end{pf*}

\begin{pf*}{Proof of Lemma~\ref{lemmal00}}
Part (a) is a well-known approximation result for tensor product
splines, see Theorem 12.7 of \cite{Schumaker2007} or Lemma 2.1 of \cite{deJonge2012}, for example. Part (b) is a direct multivariate
generalization of Lemma 1, part (b) of \cite{Shen2012}.

For part (c), note that by part (b) we have $\bolds{\theta} \in(0,1)^{J}$
such that
\[
\Biggl\|f(y|x_1,\ldots,x_d) - \sum
_{j_0=1}^{J_0}\cdots\sum_{j_{d}=1}^{J_0}
\theta_{j_0,\ldots,j_d} B_{j_0}(y) \prod_{k=1}^{d}
B_{j_k}(x_k) \Biggr\| _{\infty} \leq C_1
J^{-\beta/(d+1)}
\]
for constant $C_1 = C \|f^{(\beta)}\|_{\infty}$. Define $\bolds{\xi}$ as
the column vector of $\xi_{j_0,\ldots,j_d} =\break   \theta_{j_0,\ldots,j_d}
\int_0^1 B_{j_0}(y)\, \mathrm{d}y $ and $\mathbf{B}^{*}$ as the column vector of
$B_{\mathbf{j}}^{*}(y,\mathbf{x}) = \bar{B}_{j_0}(y) \prod_{k=1}^{d}
B_{j_k}(x_k)$. Then
%
\begin{equation}
\label{intermediateapproximation}
\bigl\|f(y|x_1,\ldots,x_d) - \bolds{
\xi}^T \mathbf{B}^{*}(y,\mathbf{x})\bigr\|_{\infty} \leq
C_1 J^{-\beta/(d+1)}.
\end{equation}
In particular, since $\|f\|_\infty<\infty$, it follows that $\|\bolds{\xi
}^T B^*\|_\infty$ is uniformly bounded.

By integration, and using the fact that B-splines add to 1, it follows that
\begin{eqnarray*}
&& \Biggl\|\sum_{j_1=1}^{J_0}\cdots\sum
_{j_{d}=1}^{J_0} \Biggl(1- \sum
_{j_0=1}^{J_0}\xi_{j_0,\ldots,j_d}\Biggr) \prod
_{k=1}^{d} B_{j_k}(x_k)\Biggr\|
_{\infty}
\\
&& \quad = \Biggl\|1 - \sum_{j_0=1}^{J_0}\cdots\sum
_{j_{d}=1}^{J_0} \xi_{j_0,\ldots,j_d} \prod
_{k=1}^{d} B_{j_k}(x_k)
\Biggr\|_{\infty}
  \leq C_1 J^{-\beta/(d+1)}
\end{eqnarray*}
for any $\mathbf{x} \in(0,1)^d$.
Applying a multivariate analog of Theorem~4.38 of \cite{Schumaker2007}
for tensor product of B-splines, we can bound the maximum norm of
coefficients in a tensor product B-spline expansion by a constant
multiple of the supremum norm of the function formed by corresponding
linear combination. This is possible by forming a dual basis consisting
of tensor product of functions in a dual basis for univariate B-splines
and by noting that the supremum norms of the elements of the dual basis
can be taken to be uniformly bounded (see Theorem~4.41 of \cite{Schumaker2007}). This leads to the relation
%
\begin{equation}
\label{approximatenormalization} \Biggl| 1- \sum_{j_0=1}^{J_0}
\xi_{j_0,\ldots,j_d} \Biggr| \leq C_1' J^{-\beta/(d+1)}
\end{equation}
for any $(j_1,\ldots,j_d) \in\{1,\ldots,J_0\}^d$ and some constant $C_1'>0$.

Define $\bolds{\eta}$ by the relations $\eta_{j_0,\ldots,j_d} = \xi
_{j_0,\ldots,j_d} / \sum_{m=1}^{J_0} \xi_{m,j_1,\ldots,j_d}$. Thus $\bolds
{\eta}\in\Delta_{J_0}^{J_0^d}$.
Then using \eqref{approximatenormalization} and the boundedness of $\|
\bolds{\xi}^T B^*\|_\infty$, we obtain
\begin{eqnarray*}
&& \bigl\|\bolds{\xi}^T B^* - \bolds{\eta}^T B^*
\bigr\|_\infty
\\
&&\quad = \sup_{\mathbf{x},y} \Biggl\llvert \sum
_{j_0=1}^{J_0}\cdots\sum_{j_{d}=1}^{J_0}
\xi_{j_0,\ldots,j_d} B_{j_0}(y) \prod_{k=1}^{d}
B_{j_k}(x_k) \Biggl[ \Biggl(\,\sum
_{m=1}^{J_0} \xi_{m,j_1,\ldots,j_d} \Biggr)^{-1}-1
\Biggr]\Biggr\rrvert
\\
&& \quad \lesssim\max_{j_1,\ldots,j_d} \Biggl|1-\sum_{m=1}^{J_0}
\xi_{m,j_1,\ldots
,j_d}\Biggr| \bigl\|\bolds{\xi}^T B^*\bigr\|_\infty
\\
&& \quad \leq C_2 J^{-\beta/(d+1)}
\end{eqnarray*}
for some positive constant $C_2$. Combining with \eqref{intermediateapproximation}, the result now follows.
\end{pf*}

\begin{pf*}{Proof of Lemma~\ref{lemmal2}}
Part (a) is a well-known approximation result for tensor Sobolev space,
see Theorem~12.7 of \cite{Schumaker2007}, for example. The proof of (b)
and (c) proceed exactly as in Lemma~\ref{lemmal00}.
\end{pf*}

\section*{Acknowledgements}
The authors would like to thank the Associate Editor and two referees
for their helpful comments that greatly improves the quality of the paper.





\printhistory
\end{document}